\documentclass[sn-mathphys]{sn-jnl}

\usepackage{subcaption}

\jyear{2022}%

\newcommand{\eps}{\varepsilon}
\newcommand{\Oes}{\Omega_\eps^s}

\newcommand{\Oest}{\Oes(t)}

\newcommand{\ue}{\mathbf{u}_\eps}
\newcommand{\hce}{\widehat{c}_\eps}
\newcommand{\hve}{\widehat{\mathbf{v}}_\eps}
\newcommand{\hD}{\widehat{\mathbf{D}}}
\newcommand{\Se}{\mathbf{S}_\eps}
\newcommand{\tx}{(t,x)}
\newcommand{\txex}{\left(t,x,\frac{x}{\eps}\right)}
\newcommand{\hx}{\widehat{x}}
\newcommand{\thx}{(t,\widehat{x})}
\newcommand{\Fe}{\mathbf{F}_\eps}
\newcommand{\Je}{J_\eps}
\newcommand{\delt}{\frac{\partial}{\partial t}}

\newcommand{\ceps}{c_\eps}
\newcommand{\tDe}{\mathbf{D}_\eps}
\newcommand{\txy}{(t,x,y)}
\renewcommand{\R}{\mathbb{R}}
\newcommand{\N}{\mathbb{N}}
\newcommand{\Ahom}{\mathbf{A}^*}
\newcommand{\Jhom}{J^*}
\newcommand{\Dhom}{\mathbf{D}^*}
\newcommand{\divy}{\nabla_y \cdot}

\newcommand{\grady}{\nabla_y}
\newcommand{\gradx}{\nabla_x}

\numberwithin{equation}{section}

\theoremstyle{thmstyleone}%
\theoremstyle{thmstyletwo}%
\newtheorem{remark}{Remark}%

\theoremstyle{thmstylethree}%

\begin{document}

\title[Transport Processes in an Elastically Deformable Perforated Medium]{Multi-Scale Modeling and Simulation of Transport Processes in an Elastically Deformable Perforated Medium}


\author*[1,2]{\fnm{Jonas} \sur{Knoch}}\email{jonas.knoch@fau.de}

\author[1]{\fnm{Markus} \sur{Gahn}}\email{markus.gahn@iwr.uni-heidelberg.de}

\author[1,2]{\fnm{Maria} \sur{Neuss-Radu}}\email{maria.neuss-radu@math.fau.de}

\author[2]{\fnm{Nicolas} \sur{Neu\ss}}\email{neuss@math.fau.de}

\affil[1]{\orgdiv{IWR}, \orgname{Ruprecht-Karls-Universit\"at Heidelberg}, \orgaddress{\street{Im Neuenheimer Feld 205}, \city{69120 Heidelberg}, \country{Germany}}}

\affil[2]{\orgdiv{Department Mathematik}, \orgname{Friedrich-Alexander-Universit\"at Erlangen-Nürnberg}, \orgaddress{\street{Cauerstra\ss e 11}, \city{91058 Erlangen}, \country{Germany}}}

\abstract{In this paper, we derive an effective model for transport processes in periodically perforated elastic media, taking into account, e.g., cyclic elastic deformations as they occur in lung tissue due to respiratory movement. The underlying microscopic problem couples the deformation of the domain with a diffusion process within a mixed \textit{Lagrangian}/\textit{Eulerian} formulation. After a transformation of the diffusion problem onto the fixed domain, we use the formal method of two-scale asymptotic expansion to derive the upscaled model, which is nonlinearly coupled through effective coefficients. The effective model is implemented and validated using an application-inspired model problem. Numerical solutions for both, cell problems and macroscopic equations, are investigated and interpreted. We use simulations to qualitatively determine the effect of the deformation on the transport process.}




\maketitle

\section{Introduction}\label{section:introduction}
Heat and mass transport in materials characterized by a complex microgeometry is an actively researched topic in various fields such as material sciences or geosciences as well as in the biological/medical field. Here, mathematical models can help to identify underlying mechanisms as well as complement or replace experiments by numerical simulations. This is of particular importance in the medical field, where experiments and animal studies or clinical trials are expensive if at all feasible. The main question of this work is, how transport processes are influenced by a complex microstructure of the carrier medium and its deformation. Applications are especially found in the context of diseases whose clinical picture includes a dysfunction of the lung, such as COVID-19, pneumonia or sepsis. Experimental evidence already indicates that in addition to the microstructure formed at the cellular level by pneumocytes, the periodic deformation caused by inhalation and exhalation also plays an important role in transport processes such as those for nutrients or respiratory gases, see e.g. \cite{Huh2010}.\\
Unfortunately, the heterogeneous microstructure of such problems can only be resolved to a very limited extent for small domains in numerical simulations. For problems involving multiple scales, one has to rely on mathematical tools that are capable of deriving approximate problems, which are easier to solve but still incorporate the information on the microgeometry. To this end, we derive in this paper an effective model for transport processes in an elasically deformable, periodically perforated medium. We approach the problem by formulating first the microscopic problem on the perforated domain in a mixed \textit{Lagrangian}/\textit{Eulerian} framework. The deformation process is described by a linear elasticity equation on the fixed microscopic domain as opposed to the transport problem, which is initially formulated on the current, moving microscopic domain. In this form, the problem is not accessible for standard upscaling tools. Therefore, we use a transformation, which is defined by means of the solution to the elasticity problem, to reformulate the transport problem on the reference domain for the elasticity equation. Then we can exploit the formal method of two-scale asymptotic expansion in order to obtain a macroscopic model describing transport processes and deformation for a periodically perforated domain. The resulting model is formulated on a homogeneous domain and involves effective coefficient functions. These effective quantities are computed by means of solutions to so-called cell problems which are 
formulated on the reference cell and 
carry the information on the microscopic geometry. Additionally, we develop and study numerical methods for the derived effective model. A lot of work has been done concerning the upscaling of fluid-structure-interaction-models, pure diffusion or elasticity, or (linear) coupled models for heat transport and deformation. In contrast, there are hardly any papers considering modeling and numerical simulations for an elasticity-diffusion problem, where the diffusion process is formulated within the \textit{Eulerian} framework and transformed into a common framework, leading to coupling through nonlinear coefficients.
We discretize both, the cell problems as well as the macroscopic problem using the finite element method and the implicit \textit{Euler} method for discretization in time. Afterwards, we use a application-inspired model problem with oscillating Dirichlet boundary conditions in the deformation  for numerical convergence tests to validate the code. For the model problem, we investigate and interpret the numerical results. Finally, again within the scope of the model problem, we give a quantitative answer to the question how the deformation affects the transport of a diffusing substance. To this end, we perform a parameter study for key parameters involved in the deformation of the domain and analyze the sensitivity of the effective model with respect to these parameters. This is accompanied by the visualization of the homogenized coefficient functions.\\
Homogenization techniques for elastic heterogeneous media have been extensively studied in \cite{Oleinik1992}. In \cite{CollisEtal2017}, a macroscopic model for transport processes, involving fluid-structure-interaction and growth has been derived by formal upscaling. However, the explicit dependence of the effective coefficients is not specified. The authors of \cite{Allaire2017a} consider the transport of electrolytes through the pore space of a deformable, charged porous medium and rigorously derive a macroscopic model, however they assume small deformations and therefore consider the transport in Lagrangian coordinates. Rigorous homogenization for reaction-diffusion models involving an evolving microstructure have been studied for example in \cite{GahnNRP,Peter07,Peter09}. In all three cases, the evolution of the microstructure is assumed to be known. For the homogenization of a model of thermoelasticity, including a given evolution of the microstructure due to phase transformation, see \cite{EdenMuntean17}. In the context of crystal precipitation and dissolution, where the evolution of the microstructure is dictated by the local dissolution/precipitation rate, an upscaled model has been derived in \cite{Noorden2009}, and the macroscopic model for a similar micro-model is analyzed in \cite{Schulzetal2016} for specific geometries leading to a dependence of the effective diffusion on the porosity. We emphasize that rigorous homogenization results for such kind of problems including a free boundary a quite rare. See for example \cite{GahnPop} for a microstructure including spherical grains, where the radii depend on the solute concentration at the surface.  \\
The paper is organized as follows: In section \ref{section:the_mathematical_model}, we introduce the microscopic problem and rewrite it in non-dimensionalized from on a fixed, heterogeneous reference domain. This is the starting point for the derivation of the effective model, which we obtain by means of a two-scale asymptotic expansion. In section \ref{section:Investigations_Concerning_the_Numerical_Simulation_of_the_Upscaled_Model}, we evaluate and interpret the numerical solutions to the macroscopic problem. Before we finish with some final remarks and an outlook, we use computer experiments in section \ref{section:sensitivity_of_transport_process} to quantitatively investigate how the cyclic deformation of the domain affects the transport process in the macroscopic model.
\section{The Mathematical Model}
\label{section:the_mathematical_model}
Let us consider a macroscopic domain $\Omega$ given by the hyper-rectangle $\Omega = (a,b)\subset \R^n$, $n\geq 2$, with $a,b \in \mathbb{Z}^n$ such that $a_i < b_i$ for all $i = 1, \dots, n$.  
This macroscopic domain contains a subdomain $\Oes$ which exhibits a periodic microstructure and can be interpreted as the solid part of a porous medium. The small parameter $\epsilon  $ with $\epsilon^{-1} \in \N$ describes the ratio between the size of the whole domain $\Omega$ and the pore size. The microscopic domain $\Oes$ is constructed in the following way: Let $Y = (0,1)^n$ be the $n$-dimensional open unit cube and $Y^s \subsetneq Y$ be a connected domain   in $Y$ that intersects all faces of $Y$, i.e. we have
$$
\partial Y^s \cap \{y_i = 0\} \neq \emptyset \quad \text{and} \quad \partial Y^s \cap \{y_i = 1\} \neq \emptyset \quad \text{ for } i = 1,...,n.
$$
Additionally, we require that opposite faces of $Y^s$ are equal, i.e. it holds that
$$
(\partial Y^s \cap \{y_i = 0 \}) + \mathbf{e}_i = (\partial Y^s \cap \{y_i = 1 \}) \quad \text{ for } i = 1,...,n.
$$
By $\Gamma := \text{int}(\partial Y^s \setminus \partial Y)$ we denote the internal boundary of $Y^s$ in $Y$.
Then we define the subdomain $\Oes$ as
$$
\Oes := \text{int}\left(\Omega \cap \left( \underset{\mathbf{k} \in \mathbb{Z}^n}{\bigcup} \eps(\overline{Y^s} + \mathbf{k}) \right) \right)
$$
and the internal boundary of the heterogeneous subdomain in $\Omega$
$$
\Gamma_\eps := \text{int}(\partial \Oes \setminus \partial \Omega).
$$
We note that by our assumptions $\Oes$ is connected. This allows us to make sense of elastic deformation and diffusion processes in $\Oes$ for $n \geq 2$.
\subsection{The Microscopic Problem}
\label{subsection:the_microscopic_problem}
Due to elastic deformation, the physical shapes of $\Oes$ and $\Omega$ change with time. We denote the current solid domain at time $t \in (0,T)$,  $T > 0$, by $\Oest$ and the current internal interface by $\Gamma_\eps(t):= \partial \Oest \setminus \partial \Omega(t)$. In the following, quantities and operators defined on or associated with the current domain $\Oest$ are marked by a hat $\widehat{\cdot}$ whereas quantities and operators defined on or associated with the fixed domain $\Oes$ do not carry a distinct label. The outer boundary of $\Oes$ is divided into a Dirichlet and a Neumann boundary for the elasticity problem, i.e. we have $\partial \Oes \cap \partial \Omega = \Gamma^{\text{elast}}_{\eps,D} \cup \Gamma^{\text{elast}}_{\eps,N}$ with $\Gamma^{\text{elast}}_{\eps,D} \cap \Gamma^{\text{elast}}_{\eps,N} = \emptyset$. Analogously, for the diffusion problem, we have $\partial \Oest \cap \partial \Omega(t) = \Gamma^{\text{diff}}_{\eps,D}(t) \cup \Gamma^{\text{diff}}_{\eps,N(t)}$ with $\Gamma^{\text{elast}}_{\eps,D}(t) \cap \Gamma^{\text{elast}}_{\eps,N}(t) = \emptyset$ for all $t \in (0,T)$.
We assume that the displacement $\ue$ and the concentration $\hce$ satisfy the system
\begin{subequations}
\begin{align}
\rho_s \frac{\partial^2 \ue}{\partial t^2} - \nabla \cdot \left( \mathbf{A} e(\ue) \right) & = \mathbf{f}^\text{elast} & \mbox{ in } & (0,T) \times \Oes,  \label{system:elasticity-diffusion-system_on_mixed_domains:a} 
\\
-  \mathbf{A} e(\ue)  \cdot \mathbf{n} & = 0 & \mbox{ on } & (0,T) \times \Gamma_\eps \cup \Gamma^{\text{elast}}_{\eps,N},
\\
\ue &= \mathbf{h} &\mbox{ on }& (0,T)\times \Gamma^{\text{elast}}_{\eps,D},
\\
\ue(0,x) = \frac{\partial \ue }{\partial t}(0,x) &= 0 &\mbox{ in }& \Oes, \label{system:elasticity-diffusion-system_on_mixed_domains:d} 
\\ 
\frac{\partial \hce}{\partial t} + \widehat{\nabla} \cdot \left( \hce \hve - \hD \widehat{\nabla} \hce \right) & = \widehat{f}^\text{diff} & \mbox{ in } & Q_\eps^T := \underset{t \in (0,T)}{\bigcup} \{t\} \times \Oest, \label{system:elasticity-diffusion-system_on_mixed_domains:e} 
\\ 
- \hD \widehat{\nabla} \hce \cdot \widehat{\mathbf{n}}  & = 0 & \mbox{ on } & G_\eps^T := \underset{t \in (0,T)}{\bigcup} \{t\} \times (\Gamma_\eps(t) \cup \Gamma^\text{diff}_{\eps,N}(t)),
\\
\hce &= \widehat{g} & \mbox{ on }&  E_\eps^T := \underset{t \in (0,T)}{\bigcup} \{t\} \times \Gamma^\text{diff}_{\eps,D}(t),
\\
\hce(0,\widehat{x}) &= \widehat{c}^0(\widehat{x}) & \mbox{ in } & \Oes(0). \label{system:elasticity-diffusion-system_on_mixed_domains:h}
\end{align}
\label{system:elasticity-diffusion-system_on_mixed_domains}
\end{subequations}
%
Here, $\rho_s > 0$ is the constant density of the solid phase, $\mathbf{A}$ is a constant fourth-order stiffness tensor and by $\mathbf{f}^\text{elast}: (0,T) \times \Omega \rightarrow \mathbb{R}^n$ we denote a body force acting on the solid phase. $e(\mathbf{w}) = \frac{1}{2}(\nabla \mathbf{w} + (\nabla \mathbf{w})^T)$ denotes the symmetric gradient and on the Dirichlet boundary of the elasticity problem, the displacement is prescribed by a function $\mathbf{h}: (0,T) \times \Gamma_D^\text{elast} \rightarrow \mathbb{R}^n$ with $\mathbf{h}(0,x) = 0$. Further, $\widehat{\mathbf{D}}$ is a constant, positive-definite diffusion tensor, $\widehat{f}^\text{diff}: \underset{t \in (0,T)}{\bigcup} \{t\} \times \Omega(t) \rightarrow \mathbb{R}$ is a source/sink term and $\widehat{g}: \underset{t \in (0,T)}{\bigcup} \{t\} \times \Gamma_D^\text{diff} \rightarrow \mathbb{R}$, $\widehat{c}^0: \Omega(t) \rightarrow \mathbb{R}$ denote the Dirichlet boundary condition and the initial condition for the diffusion problem. By $\mathbf{n}$ and $\widehat{\mathbf{n}}$, we denote the outter unit normals to the fixed domain $\Oes$ and the current domain $\Oest$, respectively. $\hve$ takes into account the velocity field that is induced by the deformation of the domain. We postpone its precise definition since we need a more detailed description of the current deformed domain first. \\ 
The elasticity equations $\eqref{system:elasticity-diffusion-system_on_mixed_domains:a}$--$\eqref{system:elasticity-diffusion-system_on_mixed_domains:d}$ are defined within the \textit{Lagrangian} framework on the fixed domain $\Oes$. In contrast, the natural setting for the diffusion equation $\eqref{system:elasticity-diffusion-system_on_mixed_domains:e}$--$\eqref{system:elasticity-diffusion-system_on_mixed_domains:h}$ is the \textit{Eulerian} framework, i.e. it is initially defined on the current deformed domain $\Oest$. In the present context, the shape of $\Oest$ is solely dictated by the elasticity problem, whereas in general, other processes, e.g. growth, might also be conceivable as driver of the deformation of the domain. Due to our assumption, we can explicitly state a transformation mapping given by the deformation
\begin{equation}
\Se(t,x) := x + \ue\tx, \quad \tx \in (0,T) \times \Oes.
\label{def:transformation_mapping}
\end{equation}
Using this definition, a point $x \in \Oes$ has the coordinates $\hx = \Se\tx $ at time $t$. Consequently, the current deformed domain $\Oest$ is given by
\begin{equation*}
\Oest := \{ \hx \in \mathbb{R}^n  \mid  \hx = \Se\tx, \; x \in \Oes \}, \quad t \in (0,T),
\end{equation*}
and analogously
\begin{equation*}
\Gamma_\eps(t) := \{ \hx \in \mathbb{R}^n  \mid  \hx = \Se\tx, \; x \in \Gamma_\eps \}, \quad t \in (0,T).
\end{equation*}
Since $\mathbf{h}(0,x)=0$, we have $\ue(0,x) = 0$ and therefore $\Oes = \Oes(0)$ and $\Gamma_\eps = \Gamma_\eps(0)$, i.e. the reference domain coincides with the initial domain. Now we are equipped to describe the velocity field $\hve$ in equation $\eqref{system:elasticity-diffusion-system_on_mixed_domains:e}$ in more detail. We assume that $\hve$ comprises only the transport of $\hce$ due to the deformation of the domain, i.e. we define 
\begin{equation*}
\hve \thx := \frac{\partial \Se}{\partial t}  (t,\Se^{-1}\thx), \quad \thx \in Q_\eps^T.
\end{equation*}
An additional contribution to the velocity field that takes into account other transport processes is of course also possible but not of interest in our case. 
\subsection{Transformation of the Diffusion Equation to the Fixed Domain}
\label{subsection:Transformation_of_the_diffusion_equation_to_the_fixed_domain}
In order to make the elasticity-diffusion system accessible for standard homogenization techniques, we transform equations $\eqref{system:elasticity-diffusion-system_on_mixed_domains:e}$--$\eqref{system:elasticity-diffusion-system_on_mixed_domains:h}$ to the fixed domain $\Oes$ with the help of the transformation mapping $\Se$ defined in $\eqref{def:transformation_mapping}$. For given, sufficiently smooth, scalar- and vector-valued functions on the current deformed space-time domain $Q_\eps^T$
\begin{equation*}
\hat{\phi}: Q_\eps^T \rightarrow \mathbb{R} \quad \text{and} \quad \hat{\mathbf{\Phi}}: Q_\eps^T \rightarrow \mathbb{R}^n,
\end{equation*}
we obtain their respective counterparts on the fixed domain via
\begin{equation*}
\phi(t,x) := \widehat{\phi}(t, \Se\tx ) \quad \text{and} \quad \mathbf{\Phi} \tx := \widehat{\mathbf{\Phi}}(t, \Se \tx ) \quad \text{ for }  (t,x) \in (0,T) \times \Oes.
\end{equation*}
Additionally, we define 
\begin{equation}
\Fe \tx := \nabla \Se \tx \quad \text{and} \quad \Je \tx := \det(\Fe\tx) \quad \text{ for }  (t,x) \in (0,T) \times \Oes. \label{def:J_eps}
\end{equation}
The following well known computation rules connect differential operators and integrals from $\Oes$ and $\Oest$ via the transformation mapping $\hat{x} = \Se\tx$. Similar transformations are used e.g. in \cite{Peter07}.
\begin{align*}
\hat{\nabla} \hat{\phi} &= \Fe^{-T}\nabla \phi,\\
\frac{\partial \hat{\phi}}{\partial t} &= \frac{\partial \phi}{\partial t} - \Fe^{-T}  \nabla \phi  \cdot \frac{\partial \Se}{\partial t},\\
\frac{\partial \Je }{\partial t}  &= \Je \mathrm{tr}\left(\Fe^{-1} \frac{\partial_t \Fe}{\partial t}\right), \\ 
\int_{\Oest} \hat{\phi} d\hat{x} &= \int_{\Oes} \phi \Je dx,\\
\int_{\partial \Oest} \hat{\mathbf{\Phi}} \cdot \hat{\mathbf{n}} d\hat{\sigma} &=
\int_{\partial \Oes}\Je  \mathbf{\Phi}\cdot \Fe^{-T} \mathbf{n}  d\sigma.
\end{align*}
%
%
Accordingly, the equations governing the diffusion process on the fixed reference domain are given by 
\begin{subequations}
\begin{align}
\delt \left( \Je \ceps \right) - \nabla \cdot \left( \tDe \nabla \ceps \right) &= \Je f^\text{diff} & \text{ in } & (0,T) \times \Oes, \label{eq:diffusion_on_fixed_domain:a}\\
-\tDe \nabla \ceps \cdot \mathbf{n} &= 0 & \text{ on } & (0,T) \times (\Gamma_\eps \cup \Gamma^\text{diff}_{\eps,N}(0)), \\
\ceps &= g & \mbox{ on } & (0,T) \times  \Gamma^\text{diff}_{\eps,D}(0)),\\
\ceps(0,x) &= c^0(x) & \mbox{ in } &\Oes(0),\label{eq:diffusion_on_fixed_domain:d}
\end{align}
\label{eq:diffusion_on_fixed_domain}
\end{subequations}
where we use the definitions $\ceps \tx := \hce(t,\Se\tx))$, $f^{\text{diff}}\tx := \widehat{f}^{\text{diff}}(t, \Se \tx))$ and
\begin{equation}
\tDe\tx := [\Je \Fe^{-1} \hD \Fe^{-T}](t,x). \label{def:D_eps}
\end{equation}
We note that the velocity field $\hve$ on the current deformed domain in equation $\eqref{system:elasticity-diffusion-system_on_mixed_domains:e}$ vanishes in the formulation on the fixed domain as it was solely generated by the movement of the domain. The information on the deformation of the domain is now encoded in the coefficients $\Je$ and $\tDe$ which depend on the displacement $\ue$ in a nonlinear manner.
\subsection{The Non-dimensional Model}
\label{subsection:dimensional_analysis_and_non_dimensionalisation}
We analyse system $\eqref{system:elasticity-diffusion-system_on_mixed_domains:a}$--$\eqref{system:elasticity-diffusion-system_on_mixed_domains:e}$, $\eqref{eq:diffusion_on_fixed_domain:a}$--$\eqref{eq:diffusion_on_fixed_domain:d}$ further by introducing dimensionless quantities. The relationship between these non-dimensional quantities gives an idea of the qualitative behaviour of the system and also might influence significantly the homogenization result. As stated in the introduction, we have in mind the experimental setup from \cite{Huh2010}, where the cross-sectional size of the microchannels is $4\cdot 10^{-4} m \times 7\cdot 10^{-5}m$. To mimic physiological breathing motion, the domain was cyclically stretched on two opposing boundaries by approximately $10\%$ to reach a maximal length of $4.4 \cdot 10^{-4}m$. In order to maintain applicability of the model to a broader variety of settings we choose as characteristic size of the domain $\bar{x} = 5 \cdot 10^{-3}$ and the characteristic cell size $l = 5 \cdot 10^{-6}m$. This results in the characteristic size of heterogeneities $\eps = l/\bar{x} = 10^{-3}$.\\
Next, let us replace the dependent and independent variables by scaled versions:
\begin{gather*}
x = x^\dagger\cdot \bar{x}, \; t = t^\dagger \cdot \bar{t}, \; \mathbf{u}_\eps = \mathbf{u}^\dagger_\eps \cdot \bar{u}, \; \mathbf{A} = \mathbf{A}^\dagger \cdot \bar{A}, \;
\mathbf{f}^\text{elast} = {\mathbf{f}^\text{elast}}^\dagger \cdot \bar{f}_e, \\
\mathbf{h} = \mathbf{h}^\dagger \cdot \bar{h}, \; \ceps = \ceps^\dagger \cdot \bar{c}, \; \tDe = J_\eps \Fe^{-1} \mathbf{D}^\dagger \cdot \bar{D} \Fe^{-T} =\tDe^\dagger \cdot \bar{D}, \\
 f^\text{diff} = {f^\text{diff}}^\dagger \cdot \bar{f}_d, \; g = g^\dagger \cdot \bar{c}, \; c^0 = {c^0}^\dagger \cdot \bar{c}.
 \label{eq:independent_variables_and_reference_quantities}
\end{gather*}
Here, quantities marked by a dagger $\cdot^\dagger$ are dimensionless variables and they carry a characteristic reference value which is marked by a bar $\bar{\cdot}$. Note that, in the definition of $\tDe$, the quantities $\Je$ and $\Fe$ do not carry a  and therefore it is sufficient to associate the non-dimensional diffusion tensor $\mathbf{D}^\dagger$ with a reference value $\bar{D}$. Then, with the definition $\tDe^\dagger := \Je \Fe^{-1} \mathbf{D}^\dagger  \Fe^{-T}$, we can write $\tDe = \tDe^\dagger \cdot \bar{D}$. The introduction of the non-dimensional spatial variables implies that the equations are now posed on the non-dimensional domain ${\Omega_\eps^s}^\dagger$ of unit size.
Realistic values for the reference quantities according to the application we have in mind can be found in Table \ref{tab:non-dimensionalisation_reference_quantities}. 
\begin{table}
\centering
\begin{tabular}{r  r r p{8cm} }
 quantity  & value & unit & comment \\ \hline \hline
$\eps$ & $10^{-3}$ &  - & proportion of microstructure and macro-domain \\
$\bar{x}$ & $5\cdot 10^{-3}$  & $m$  &  characteristic domain size    \\
$\bar{t}$  & $1$  &   $s$  &    characteristic time for one breath     \\
$\bar{u}$, $\bar{h}$  & $5\cdot 10^{-4}$    &  $m$   &    characteristic displacement \\
$\bar{A}$ & $1 \cdot 10^{4}$  & $Pa$   &      characteristic elastic modulus  \\
$\bar{f}_e$ & $2 \cdot 10^5$ & $\frac{N}{m^3}$ & characteristic body force density \\
$\bar{c}$, $\bar{g}$ & $1$ & $\frac{mol}{m^3}$ & characteristic concentration \\
$\bar{D}$ & $10^{-4} - 10^{-11}$ & $\frac{m^2}{s}$ & range for characteristic diffusion coefficient\\
$\bar{f}_d$ & $1$ & $\frac{mol}{m^3 s}$ & characteristic reaction rate of source/sink term \\ \hline %
\end{tabular}
\caption{{\bf Characteristic reference quantities} as they are used for the non-dimensionalisation of system $\eqref{system:elasticity-diffusion-system_on_mixed_domains:a}$--$\eqref{system:elasticity-diffusion-system_on_mixed_domains:e}$, $\eqref{eq:diffusion_on_fixed_domain:a}$--$\eqref{eq:diffusion_on_fixed_domain:d}$. Note that the value for $\bar{f}_e$ is chosen in such a way that the right-hand side of the elasticity equation is of order 1 and a non-zero force-term in the macroscopic equation occurs. The value might not be realistic for the application we have in mind but is chosen nonetheless to sustain generality of the model. For more realistic (smaller) values of $\bar{f}_e$, the force term would not occur in the macroscopic model.}
\label{tab:non-dimensionalisation_reference_quantities}
\end{table}
Using the relations in $\eqref{eq:independent_variables_and_reference_quantities}$ and dropping the daggers for the sake of an easier notation we obtain
\begin{align*}
\frac{\rho_s \bar{x}^2}{\bar{t}^2 \bar{A}} \frac{\partial^2\ue}{\partial t^2}  - \nabla \cdot \left( Ae(\ue) \right) &= \frac{\bar{f}_e \bar{x}^2}{\bar{A} \bar{u}}  \mathbf{f}^\text{elast},\\
\delt \left(\Je \ceps \right) - \frac{\bar{D} \bar{t}}{\bar{x}^2} \nabla \cdot \left( \tDe \nabla \ceps \right) &= \frac{\bar{f}_d \bar{t}}{\bar{c}} f^\text{diff}.
\end{align*}
According to Table \ref{tab:non-dimensionalisation_reference_quantities}, we compute (since $\eps = l/\bar{x} = 10^{-3}$)
\begin{gather*}
\frac{\rho_s \bar{x}^2}{\bar{t}^2 \bar{A}} \approx O(\eps^2), 
\quad \frac{\bar{f}_e \bar{x}^2}{\bar{A} \bar{u}} \approx O(1),
\quad \frac{\bar{D} \bar{t}}{\bar{x}^2} \approx O(1), 
\quad \frac{\bar{f}_d \bar{t}}{\bar{c}} \approx O(1).
\end{gather*} 
and therefore obtain the non-dimensionalised elasticity-diffusion system on the non-dimensional, fixed domain:
\begin{subequations}
\begin{align}
\eps^2 \frac{\partial^2 \ue}{\partial t^2} - \nabla \cdot \left( \mathbf{A} e(\ue) \right) & = \mathbf{f}^\text{elast} & \mbox{ in } & (0,T) \times \Oes, 
\label{system:elasticity-diffusion-system_on_fixed_domain:a}
\\
-  \mathbf{A} e(\ue)  \cdot \mathbf{n} & = 0 & \mbox{ on } & (0,T) \times (\Gamma_\eps \cup \Gamma^\text{elast}_{\eps,N}),
\\
\ue &= \mathbf{h} & \mbox{ on } & (0,T) \times \Gamma^\text{elast}_{\eps,D},
\\
\ue(0,x) = \frac{\partial \ue }{\partial t}(0,x) &= 0 &\mbox{ in }& \Oes, \label{system:elasticity-diffusion-system_on_fixed_domain:d}
\\
\delt \left( \Je \ceps \right) - \nabla \cdot \left( \tDe \nabla \ceps \right) &= \Je f^\text{diff} & \text{ in } & (0,T) \times \Oes,
\label{system:elasticity-diffusion-system_on_fixed_domain:e}
\\
-\tDe \nabla \ceps \cdot \mathbf{n} &= 0 & \text{ on } & (0,T) \times (\Gamma_\eps \cup \Gamma^\text{diff}_{\eps,N}),
\\
\ceps &= g & \mbox{ on } & (0,T) \times \Gamma^\text{diff}_{\eps,D},
\\
\ceps(0,x) &= c^0(x) & \mbox{ in } & \Oes. \label{system:elasticity-diffusion-system_on_fixed_domain:h}
\end{align}
\label{system:elasticity-diffusion-system_on_fixed_domain}
\end{subequations}

\begin{remark}
The existence of a weak solution for the elasticity equation $\eqref{system:elasticity-diffusion-system_on_fixed_domain:a}$ is quite standard. The crucial point is the existence for the reaction-diffusion equation $\eqref{system:elasticity-diffusion-system_on_fixed_domain:e}$, since the coefficients depend on the microscopic displacement $\mathbf{u}_{\eps}$ and may degenerate in space and time for  vanishing $\Je$. To overcome this problem higher regularity for the displacement is necessary, to guarantee the positivity of $\Je$ as well as the injectivity of  $\mathbf{u}_{\eps}$, at least locally in time. However, this local existence interval may depend on $\epsilon$ and could vanish for $\epsilon \to 0$. Even for finite $\epsilon$ strong assumptions and compatibility conditions on the data (and especially the microscopic domain) are necessary to obtain enough regularity for the displacement. Hence, the rigorous analytical treatment of this problem is an open question.

\end{remark}

\subsection{Upscaling in the Fixed Domain}
\label{subsection:upscaling_in_the_fixed_domain}
Although the domain for the elasticity-diffusion system is now fixed in time, we still have to deal with its heterogeneous microstructure. This issue is addressed by application of the formal two-scale asymptotic expansion method with the objective of deriving effective equations on the spatially homogeneous macroscopic domain $\Omega$ and complementary cell problems on the scaled microscopic domain $Y^s$. The cell problems account for the heterogeneous nature of the initial problem and will give rise to the definition of effective coefficient functions for the macroscopic problem.
\subsubsection{Upscaling of the Elasticity Problem}
Following the standard approach of two-scale asymptotic expansion, we postulate expansions of the form
\begin{equation}
\ue\tx = \sum_{i = 0}^\infty \eps^i \mathbf{u}_i \txex, \quad \ceps \tx = \sum_{i = 0}^\infty \eps^i c_i \txex,
\label{eq:two-scale-expansion}
\end{equation}
for the microscopic displacement $\ue$ and  concentration  $\ceps$,
with $y := \frac{x}{\eps}$, a new spatial variable that is associated with the microscopic scale of the problem. Here, each function $\mathbf{u}_i \txy$ and $ c_i \txy$ is a function of time $t \in (0,T)$, the macroscopic spatial variable $x \in \Omega$, and the microscopic spatial variable $y \in Y^s$ and is $Y$-periodic with respect to $y$. For derivatives of such functions, we have by the chain rule
$$
\nabla = \gradx + \frac{1}{\eps} \grady, \qquad e = e_x + \frac{1}{\eps} e_y,
$$
and analogous rules hold for the divergence operator. We will see that the expansion for $\ue$ implies a two-scale asymptotic expansion of type  
\begin{equation}
    \Psi_\eps \tx = \sum_{i = 0}^\infty \eps^i \Psi_i \txex, \label{eq:general_two_scale_expansion}
\end{equation}
where $\Psi_\eps$ and $\Psi_i, i=1,2,...,$ can be scalar- or vector-valued, also for the coefficients $\Je$ and $\tDe$.
We start with the derivation of a macroscopic equation for the displacement, and identify the zeroth order term $\mathbf{u}_0$ and the first order term $\mathbf{u}_1$ in the expansion $\eqref{eq:two-scale-expansion}$ for the microscopic displacement $\ue$. We emphasize that these results are well known in the literature, even for rigorous homogenization, see for example \cite{Oleinik1992}. However, for the sake of completeness and since we need the exiplicit representations for $u_1$, we give some details for the upscaling process.
Hence, let  us  plug in the expansion for $\ue$ from (\ref{eq:two-scale-expansion}) into equations $\eqref{system:elasticity-diffusion-system_on_fixed_domain:a}$--$\eqref{system:elasticity-diffusion-system_on_fixed_domain:d}$. With the computation rules above we obtain a series in $\eps$ for the elasticity equation:
\begin{align}
\begin{split}
\eps^2 \frac{\partial^2}{\partial t^2} & \left( \mathbf{u}_0 \txy + \eps\mathbf{u}_1 \txy + ... \right) \\
& - \eps^{-2} \left[ \grady \cdot \left( \mathbf{A} e_y(\mathbf{u}_0) \right) \right]\txy \\
& - \eps^{-1} \left[ \grady \cdot \left(\mathbf{A} ( e_y(\mathbf{u}_1) + e_x(\mathbf{u}_0)) \right)  + \gradx \cdot \left( \mathbf{A} e_y( \mathbf{u}_0) \right) \right] \txy \\
& - \eps^0 \left[ \grady \cdot \left(\mathbf{A} ( e_y(\mathbf{u}_2) + e_x(\mathbf{u}_1)) \right)  + \gradx \cdot \left( \mathbf{A} (e_y( \mathbf{u}_1) + e_x(\mathbf{u}_0)) \right) \right] \txy \\
& - ...\\
& = \mathbf{f}^\text{elast}\tx,
\end{split}
\label{eq:elasticty_cascade}
\end{align}
and for the boundary conditions
\begin{align*}
& - \eps^{-1} \left[\mathbf{A}e_y(\mathbf{u}_0) \cdot \mathbf{n} \right] \txy \\
& - \eps^0 \left[\mathbf{A}\left( e_y(\mathbf{u}_1)  + e_x(\mathbf{u}_0) \right) \cdot \mathbf{n} \right] \txy \\
& - \eps^1 \left[\mathbf{A}\left( e_y(\mathbf{u}_2)  + e_x(\mathbf{u}_1) \right) \cdot \mathbf{n} \right] \txy \\
& - ...  = 0, 
\end{align*}
and
\begin{equation*}
\left[\mathbf{u}_0 + \eps \mathbf{u}_1+ \eps^2 \mathbf{u}_2  + ... \right] \txy = \mathbf{h}\tx,
\end{equation*}
as well as for the initial condition 
\begin{equation*}
\left[\mathbf{u}_0 + \eps \mathbf{u}_1+ \eps^2 \mathbf{u}_2  + ... \right] (0,x,y) = \frac{\partial}{\partial t}\left[\mathbf{u}_0 + \eps \mathbf{u}_1+ \eps^2 \mathbf{u}_2  + ... \right] (0,x,y) = 0.
\end{equation*}
By comparing the coefficients in equation $\eqref{eq:elasticty_cascade}$, we obtain from the terms of order $\epsilon^{-2}$  together with the boundary condition of order $\epsilon^{-1}$ the following cell problem for $\mathbf{u}_0$:
\begin{subequations}
\begin{align*}
-\grady \cdot \left( \mathbf{A}e_y(\mathbf{u}_0) \right) &= 0 & \mbox{ in } & Y^s, \\
- \mathbf{A}e_y(\mathbf{u}_0) \cdot \mathbf{n} &= 0 & \mbox{ on } & \Gamma, \\
\mathbf{u}_0 \text{ is } Y^s & \text{-periodic in } y.
\end{align*}
\label{eq:problem_for_u0}
\end{subequations}
Since $\mathbf{A}$ is positive on the space of symmetric matrices we obtain $e_y(\mathbf{u}_0) =0$ and therefore, since $Y^s$ is connected, $\mathbf{u}_0$ is a rigid displacement with respect to $y$. However, the only periodic rigid displacements are constant functions and therefore $\mathbf{u}_0\txy = \mathbf{u}\tx$. Since $e_y(\mathbf{u}) = 0$, the term of order $\eps^{-1}$ in $\eqref{eq:elasticty_cascade}$ yields the problem
\begin{subequations}
\begin{align*}
- \grady \cdot \left[ \mathbf{A} \left( e_y(\mathbf{u}_1) + e_x(\mathbf{u}) \right) \right] &= 0 & \mbox{ in } & Y^s, \\ 
- \mathbf{A} \left( e_y(\mathbf{u}_1) + e_x(\mathbf{u}) \right) \cdot \mathbf{n} &= 0 & \mbox{ on } & \Gamma, \\
\mathbf{u}_1 \text{ is } Y^s & \text{-periodic in } y,
\end{align*}
\label{eq:problem_for_u1}
\end{subequations}
i.e. $\mathbf{u}_1$ can be determined in terms of $\mathbf{u}$ and due to the linearity of the equation we have the representation (unique up to a constant depending on $t$ and $x$, which we choose equal to zero by assuming mean value zero for $\mathbf{u}_1$)
\begin{equation}
\mathbf{u}_1\txy = \sum_{i,j = 1}^n e_x(\mathbf{u})_{ij}(t,x) \boldsymbol{\chi}_{ij}(y),
\label{eq:u1_representation}
\end{equation}
where $\boldsymbol{\chi}_{ij}$, for $i,j = 1,...,n$ are the solutions to the cell problems 
\begin{subequations}
\begin{align}
- \nabla_y \cdot \left[\mathbf{A}\left( \frac{ \mathbf{e}_i \otimes \mathbf{e}_j + \mathbf{e}_j \otimes \mathbf{e}_i}{2} + e_y(\boldsymbol{\chi}_{ij}) \right)\right] &= 0 &\mbox{ in }& Y^s, \label{eq:elasticity_cell_problems:a}
\\
- \mathbf{A}\left( \frac{ \mathbf{e}_i \otimes \mathbf{e}_j + \mathbf{e}_j \otimes \mathbf{e}_i}{2} + e_y(\boldsymbol{\chi}_{ij}) \right) \cdot \mathbf{n} &= 0 &\mbox{ on }& \Gamma,
\\
\boldsymbol{\chi}_{ij} \mbox{ is } Y \mbox{-periodic in $y$, } \, \int_{Y^s} \boldsymbol{\chi}_{ij} dy = 0.
\end{align}
\label{eq:elasticity_cell_problems}
\end{subequations}
Here, $\mathbf{e}_i$ is the $i$-th canonical unit basis vector in $\mathbb{R}^n$ and $(\mathbf{a} \otimes \mathbf{b})_{ij} = a_ib_j$ is the dyadic product of two vectors $\mathbf{a}, \mathbf{b} \in \mathbb{R}^n$.
Finally, the term of order $\eps^0$, together with boundary conditions, gives an equation for $\mathbf{u}_2$ in terms of $\mathbf{u}$ and $\mathbf{u}_1$: 
\begin{subequations}
\begin{align*}
-\grady \cdot \left(\mathbf{A} ( e_y(\mathbf{u}_2) + e_x(\mathbf{u}_1)) \right) &= \gradx \cdot \left( \mathbf{A} (e_y( \mathbf{u}_1) + e_x(\mathbf{u})) \right) + \mathbf{f}^\text{elast} & \mbox{ in } & Y^s, \\
- \mathbf{A} \left( e_y(\mathbf{u}_2) + e_x(\mathbf{u}_1) \right) \cdot \mathbf{n} &= 0 & \mbox{ on } & \Gamma, \\
\mathbf{u}_2 \text{ is } Y^s & \text{-periodic in } y.
\end{align*}
\label{eq:problem_for_u2}
\end{subequations}
Integration over $Y^s$ and integration by parts with respect to $y$ gives the macroscopic problem:
\begin{align*}
- \gradx \cdot\int_{Y^s}   \mathbf{A} (e_y( \mathbf{u}_1) +
  e_x(\mathbf{u}))  dy
  &= \lvert Y^s\rvert \mathbf{f}^\text{elast} & \mbox{ in } & (0,T) \times \Omega,
\end{align*}
which is actually an equation for $\mathbf{u}$. In fact, exploiting now the representation (\ref{eq:u1_representation}) of $\mathbf{u}_1$, we see that $\mathbf{u}$ solves the quasi-static problem
\begin{subequations}
\begin{align}
-\nabla \cdot \left( \Ahom e( \mathbf{u}) \right) &= \lvert Y^s\rvert \mathbf{f}^\text{elast} \quad &\mbox{ in }& (0,T) \times \Omega, \label{eq:homogenized_elasticity_problem:a}
\\
-\Ahom e( \mathbf{u}) \cdot \mathbf{n} &= 0 & \mbox{ in } & (0,T) \times \Gamma^\text{elast}_N,
\\
\mathbf{u} &= \mathbf{h} &\mbox{ on }& (0,T)\times \Gamma^\text{elast}_{D}, \label{eq:homogenized_elasticity_problem:c}
\end{align}
\label{eq:homogenized_elasticity_problem}
\end{subequations}
where the effective elasticity tensor $\Ahom$ is given by means of the solutions to the cell problems:
\begin{equation}
\Ahom_{ijrs} = \sum_{k,l=1}^n \int_{Y^s} \mathbf{A}_{ijkl}\left( \delta_{kr}\delta_{ls} - e_y(\boldsymbol{\chi}_{rs})_{kl} \right) dy \quad \mbox{ for } i,j,r,s = 1,...,n.
\label{Ahom}
\end{equation}
\subsubsection{Upscaling of the Diffusion Problem}
For the upscaling of the transport problem $\eqref{system:elasticity-diffusion-system_on_fixed_domain:a}$--$\eqref{system:elasticity-diffusion-system_on_fixed_domain:e}$, we first expand the coefficients  $\Je$ and $\tDe$ with respect to $\epsilon$ according to $\eqref{eq:general_two_scale_expansion}$ by using the expansion of $\ue$ and the Taylor expansion resp. Neumann series. We start with the expansion for $\Se$, which is given by
\begin{equation*}
\Se \tx = x + \mathbf{u}\tx + \eps \mathbf{u}_1 \txex + \eps^2 \mathbf{u}_2 \txex + ...,
\end{equation*}
and therefore, using also representation $\eqref{eq:u1_representation}$ for $\mathbf{u}_1$, we get
\begin{equation*}
\begin{split}
\Fe \tx = \nabla \Se \tx = \mathbf{E}_n + \gradx \mathbf{u}\tx + \sum_{i,j = 1}^n e_x(\mathbf{u})_{ij}\tx \grady \boldsymbol{\chi}_{ij}\left(\frac{x}{\eps} \right) \\
+ \eps \left( \gradx \mathbf{u}_1 + \grady \mathbf{u}_2 \right) \txex + \eps^2 \left( \gradx \mathbf{u}_2 + \grady \mathbf{u}_3 \right) \txex + ...,
\end{split}
\end{equation*}
with $\mathbf{E}_n$ the unit matrix in $ \mathbb{R}^{n\times n}$, which is also an expansion of the form $\eqref{eq:general_two_scale_expansion}$. Hence, we have  
\begin{align*}
   \mathbf{F}_0(t,x,y) = \mathbf{E}_n + \gradx \mathbf{u}\tx + \sum_{i,j = 1}^n e_x(\mathbf{u})_{ij}\tx \grady \boldsymbol{\chi}_{ij}(y) . 
\end{align*}
Since we use a formal upscaling approach, we can assume that $\mathbf{F}_0(t,x,y)$ is invertible. However, we emphasize that for a rigorous homogenization process this would be a critical point.
As the determinant and the inverse are nonlinear operations, we employ the following linearizations: For general  matrices $\mathbf{A}, \mathbf{B} \in \mathbb{R}^{n \times n}$ with $\mathbf{A}$ non-singular and  $\eps > 0$ small enough, one has with the Taylor expansion for $\det$ and the Neumann series for the inverse
\begin{align*}
\det( \mathbf{A} + \eps \mathbf{B}) &= \det(\mathbf{A}) + O(\eps),\\
( \mathbf{A} + \eps \mathbf{B})^{-1} &= \mathbf{A}^{-1} + O(\eps).
\end{align*}
Consequently, we obtain for the determinant $\Je$ and the diffusion tensor $\tDe$ 
\begin{align*}
\Je \tx &= \det( \Fe \tx ) 
= \det\left(\mathbf{F}_0\txex \right) + O(\eps),\\
\tDe \tx &= [\Je \Fe^{-1} \hD \Fe^{-T}] \tx 
= [J_0 \mathbf{F}_0^{-1} \hD \mathbf{F}_0^{-T}] \txex + O(\eps).
\end{align*}
In the following we will see that for the macroscopic model for the transport equation it will be only necessary to consider in detail the leading order terms 
%
\begin{subequations}
\begin{align}
J_0 \txy &= \det \left( \mathbf{F}_0 \txy \right)  \label{eq:J0_and_D0:a},\\
\mathbf{D}_0 \txy &= [J_0 \mathbf{F}_0^{-1} \hD \mathbf{F}_0^{-T}] \txy, \label{eq:J0_and_D0:b}
\end{align}
\label{eq:J0_and_D0}
\end{subequations}
of the expansions for $\Je$, $\tDe$, as the higher order terms will not play a role in the homogenized equation.  Since $J_0$ is bounded from below by a positive constant, we obtain that $\mathbf{D}_0$ is positive.
Now we are equipped to plug the expansions for $\ceps$ (\ref{eq:two-scale-expansion}), $\Je$ and $\tDe$ into equations $\eqref{system:elasticity-diffusion-system_on_fixed_domain:e}$--$\eqref{system:elasticity-diffusion-system_on_fixed_domain:h}$. Analogously to the asymptotic expansion for the elasticity subproblem $\eqref{system:elasticity-diffusion-system_on_fixed_domain:a}$--$\eqref{system:elasticity-diffusion-system_on_fixed_domain:d}$ before, we obtain 
\begin{align}
\begin{split}
\frac{\partial }{\partial t} &[(J_0 + \eps J_1 + \eps^2 J_2 + ...)( c_0 + \eps c_1 + \eps^2 c_2 + ...)] \txy \\
& - \eps^{-2} \left[ \grady \cdot \left( \mathbf{D}_0 \grady c_0 \right) \right] \txy \\
& - \eps^{-1} \left[ \grady \cdot \left( \mathbf{D}_0( \grady c_1 + \gradx c_0)) + \mathbf{D}_1 \grady c_0 \right) + \gradx \cdot \left( \mathbf{D}_0 \grady c_0 \right) \right] \txy \\
& - \eps^{0} [ \grady \cdot \left( \mathbf{D}_0( \grady c_2 + \gradx c_1) + \mathbf{D}_1( \grady c_1 + \gradx c_0) + \mathbf{D}_2 \grady c_0 \right) + \\
& \hspace{1cm}\gradx \cdot \left( \mathbf{D}_0( \grady c_1 + \gradx c_0) + \mathbf{D}_1 \grady c_0 \right)] \txy \\
& - ... \\
& = (J_0 + \eps J_1 + \eps^2 J_2 + ...) f^\text{diff}, \tx
\end{split}
\label{eq:diffusion_cascade} 
\end{align}
together with the boundary conditions
\begin{align*}
\begin{split}
& - \eps^{-1} \left[ \mathbf{D}_0 \grady c_0 \cdot \mathbf{n} \right] \txy \\
& - \eps^0 \left[ \left( \mathbf{D}_0 ( \grady c_1 + \gradx c_0) + \mathbf{D}_1 \grady c_0 \right) \cdot \mathbf{n} \right] \txy \\
& - \eps^1 \left[ \left( \mathbf{D}_0( \grady c_2 + \gradx c_1) + \mathbf{D}_1( \grady c_1 + \gradx c_0) + \mathbf{D}_2 \grady c_0 \right) \cdot \mathbf{n} \right] \txy\\
& - ... = 0,
\end{split}
\end{align*}
and
\begin{equation*}
[c_0 + \eps c_1 + \eps^2 c_2 + ...]\txy = g\tx,
\end{equation*}
as well as the initial condition 
\begin{equation*}
[c_0 + \eps c_1 + \eps^2 c_2 + ...](0,x,y) = c^0 (x).
\end{equation*}
In the same spirit as for the elasticity problem, we obtain for order $\eps^{-2}$ from equation $\eqref{eq:diffusion_cascade}$ the problem
\begin{subequations}
\begin{align*}
\grady \cdot \left( \mathbf{D}_0 \grady c_0 \right) &= 0 & \mbox{ in } & (0,T) \times \Omega \times Y^s, \\
- \mathbf{D}_0 \grady c_0 \cdot \mathbf{n} &= 0 & \mbox{ on } & (0,T) \times \Omega \times \Gamma, \\
c_0 \text{ is } Y^s & \text{-periodic in } y.
\end{align*}
\end{subequations}
The positivity of $\mathbf{D}_0$ implies that $c_0 \txy = c\tx $.
Therefore the problem for $c_{\epsilon}$ of order $\eps^{-1}$ reduces to
\begin{subequations}
\begin{align*}
- \grady \cdot \left[ \mathbf{D}_0 ( \grady c_1 + \gradx c ) \right] & = 0 & \mbox{ in } & (0,T) \times \Omega \times Y^s, \\
- \mathbf{D}_0 (\grady c_1 + \gradx c) \cdot \mathbf{n} & = 0 & \mbox{ on } & (0,T) \times \Omega \times \Gamma, \\
c_1 \text{ is } Y^s & \text{-periodic in } y.
\end{align*}
\end{subequations}
Linearity of this equations implies again that we can write down a representation of $c_1$ in terms of the gradient of $c$ and the solution of diffusion cell problems $\eta_i$, 
\begin{equation}
c_1 \txy = \sum_{i = 1}^n \partial_i c(t,x) \eta_i(t,x,y),
\label{eq:c1_representation}
\end{equation}
where the $n$ cell problems for $\eta_i$ are given by
\begin{subequations}
\begin{align}
-\divy \left[\mathbf{D}_0\txy(\mathbf{e}_i + \grady \eta_i\txy)\right] &= 0 &\mbox{ in }& (0,T) \times \Omega \times Y^s, \label{eq:diffusion_cell_problems:a} \\
- \mathbf{D}_0\txy[\mathbf{e}_i + \grady \eta_i\txy] \cdot \mathbf{n} &= 0 &\mbox{ on } & (0,T) \times \Omega \times \Gamma, \\
\eta_i \text{ is $Y^s$-periodic in } y, \quad \int_{Y^s} \eta_i\txy dy = 0, \label{eq:diffusion_cell_problems:c}
\end{align}
\label{eq:diffusion_cell_problems}
\end{subequations}
for $i = 1,...,n$. Eventually, we obtain an equation for $c_2$ in terms of $c$ and $c_1$ from the term of order $\eps^0$ in $\eqref{eq:diffusion_cascade}$: 
\begin{align*}
- \grady \cdot \left[ \mathbf{D}_0  ( \grady c_2 + \gradx c_1) + \mathbf{D}_1( \grady c_1 + \gradx c)\right] \textcolor{white}{------} \\ = - \frac{\partial}{\partial t} (J_0 c) + \gradx \cdot \left[ \mathbf{D}_0(\grady c_1 + \gradx c) \right]  + J_0 f^\text{diff}  &
& \mbox{ in }  & (0,T) \times \Omega \times Y^s, \\
-  \left( \mathbf{D}_0  ( \grady c_2 + \gradx c_1) + \mathbf{D}_1( \grady c_1 + \gradx c) \right) \cdot \mathbf{n}  & = 0  &\mbox{ on } & (0,T) \times \Omega \times \Gamma, \\
c_2 \text{ is } Y^s \text{-periodic in } & y.
\end{align*}
Integrating again over $Y^s$ we obtain the compatibility condition
\begin{align*}
\frac{\partial}{\partial t} \left(\int_{Y^s}J_0 dy \; c \right) - \gradx \cdot \int_{Y^s} \mathbf{D}_0(\grady c_1 + \gradx c)dy & = \int_{Y^s}J_0 dy \; f^\text{diff} & \mbox{ in } (0,T) \times \Omega.
\end{align*}
Let us define the effective coefficients $\Jhom$ and $\Dhom$ via
\begin{subequations}
\begin{align}
\Jhom\tx &= \int_{Y^s} J_0\txy dy,\\
\Dhom_{ij}\tx &= \sum_{k = 1}^n \int_{Y^s} \mathbf{D}_{0,ik}\txy \left( \delta_{kj} + \frac{\partial}{\partial y_k} \eta_j\txy \right) dy \quad \text{ for } i,j = 1,...,n.
\end{align}
\label{eq:Jhom_Dhom_definition}
\end{subequations}
By using representation $\eqref{eq:c1_representation}$ for $c_1$ we obtain after an elemental calculation that $c$ solves the following macroscopic problem:
\begin{subequations}
\begin{align}
\frac{\partial}{\partial t} \left(\Jhom c \right) - \nabla \cdot (\Dhom \nabla c) & = \Jhom f^\text{diff} & \mbox{ in } & (0,T) \times \Omega, \label{homogenized_diffusion_problem:a} \\
- \Dhom \nabla c \cdot \mathbf{n} &= 0 & \mbox{ on } & (0,T) \times \Gamma^\text{diff}_N,\\
c\tx &= g & \mbox{ on } & (0,T) \times \Gamma^\text{diff}_D,\\
c(0,x)& = c^0(x) & \mbox{ in } & \Omega. \label{homogenized_diffusion_problem:d}
\end{align}
\label{eq:homogenized_diffusion_problem}
\end{subequations}
Note that, even under the assumption that $\widehat{\mathbf{D}}$ in the initial problem (\ref{system:elasticity-diffusion-system_on_mixed_domains}e) is constant, we obtain now time- and space-dependent homogenized coefficients $\Jhom, \Dhom$ for the diffusion problem, since they carry the information on the deformation of the domain after the transformation into the fixed setting. 
\subsection{The Upscaled Equations on the Fixed Domain}
\label{subsection:The_Upscaled_Equations_on_the_Fixed_Domain}
We now summarize the upscaled equations which we have derived previously. They consist of coupled, effective, macroscopic problems for the displacement and the concentration, and are complemented by cells problems for the elasticity problem and for the diffusion problem, respectively. 
Consider the macroscopic problem
\begin{subequations}
\begin{align*}
-\nabla \cdot \left( \Ahom e( \mathbf{u}) \right) &= \lvert Y^s\rvert \mathbf{f}^\text{elast} \quad &\mbox{ in }& (0,T) \times \Omega,
\\
-\Ahom e( \mathbf{u}) \cdot \mathbf{n} &= 0 & \mbox{ in } & (0,T) \times \Gamma^\text{elast}_N,
\\
\mathbf{u} &= \mathbf{h} &\mbox{ on }& (0,T)\times \Gamma^\text{elast}_{D}, \\
\frac{\partial}{\partial t} \left(\Jhom c \right) - \nabla \cdot (\Dhom \nabla c) & = \Jhom f^\text{diff} & \mbox{ in } & (0,T) \times \Omega, \\
- \Dhom \nabla c \cdot \mathbf{n} &= 0 & \mbox{ on } & (0,T) \times \Gamma^\text{diff}_N,\\
c&= g & \mbox{ on } & (0,T) \times \Gamma^\text{diff}_D,\\
c(0)& = c^0(x) & \mbox{ in } & \Omega.
\end{align*}
\end{subequations}
The elasticity cell solutions $\boldsymbol{\chi}_{ij}$, for $i,j = 1,...,n$ are obtained by solving 
\begin{subequations}
\begin{align*}
- \nabla_y \cdot \left[\mathbf{A}\left( \frac{ \mathbf{e}_i \otimes \mathbf{e}_j + \mathbf{e}_j \otimes \mathbf{e}_i}{2} + e_y(\boldsymbol{\chi}_{ij}(y)) \right)\right] &= 0 &\mbox{ in }& Y^s,
\\
- \mathbf{A}\left( \frac{ \mathbf{e}_i \otimes \mathbf{e}_j + \mathbf{e}_j \otimes \mathbf{e}_i}{2} + e_y(\boldsymbol{\chi}_{ij}(y)) \right) \cdot \mathbf{n} &= 0 &\mbox{ on }& \Gamma,
\\
\boldsymbol{\chi}_{ij} \mbox{ is } Y^s \mbox{-periodic in $y$, } \, \int_{Y^s} \boldsymbol{\chi}_{ij}(y) dy = 0,
\end{align*}
\end{subequations}
and the solutions to the diffusion cell problems $\eta_i$, $i = 1,...,n$ are obtained from
\begin{subequations}
\begin{align*}
-\divy \left[\mathbf{D}_0\txy(\mathbf{e}_i + \grady \eta_i\txy)\right] &= 0 &\mbox{ in }& (0,T) \times \Omega \times Y^s, \\
- \mathbf{D}_0\txy[\mathbf{e}_i + \grady \eta_i\txy] \cdot \mathbf{n} &= 0 &\mbox{ on } & (0,T) \times \Omega \times \Gamma ,\\
\eta_i \text{ is $Y^s$-periodic in } y, \quad \int_{Y^s} \eta_i\txy dy = 0.
\end{align*}
\end{subequations}
From the cell problems, the following effective quantities are computed:
\begin{subequations}
\begin{align*}
\Ahom_{ijrs} &= \sum_{k,l=1}^n \int_{Y^s} \mathbf{A}_{ijkl}\left( \delta_{kr}\delta_{ls} - e_y(\boldsymbol{\chi}_{rs}(y))_{kl} \right) dy \quad \mbox{ for } i,j,r,s = 1,...,n, \\
\Jhom\tx &= \int_{Y^s} J_0\txy dy,\\
\Dhom_{ij}\tx &= \sum_{k = 1}^n \int_{Y^s} \mathbf{D}_{0,ik}\txy \left( \delta_{kj} + \frac{\partial}{\partial y_k} \eta_j\txy \right) dy \quad \text{ for } i,j = 1,...,n.
\end{align*}
\end{subequations}
\section{Numerical Investigations of the Upscaled Model}
\label{section:Investigations_Concerning_the_Numerical_Simulation_of_the_Upscaled_Model}
\begin{figure}
\centering
\begin{subfigure}[t]{0.49\textwidth}
\includegraphics[width = \textwidth]{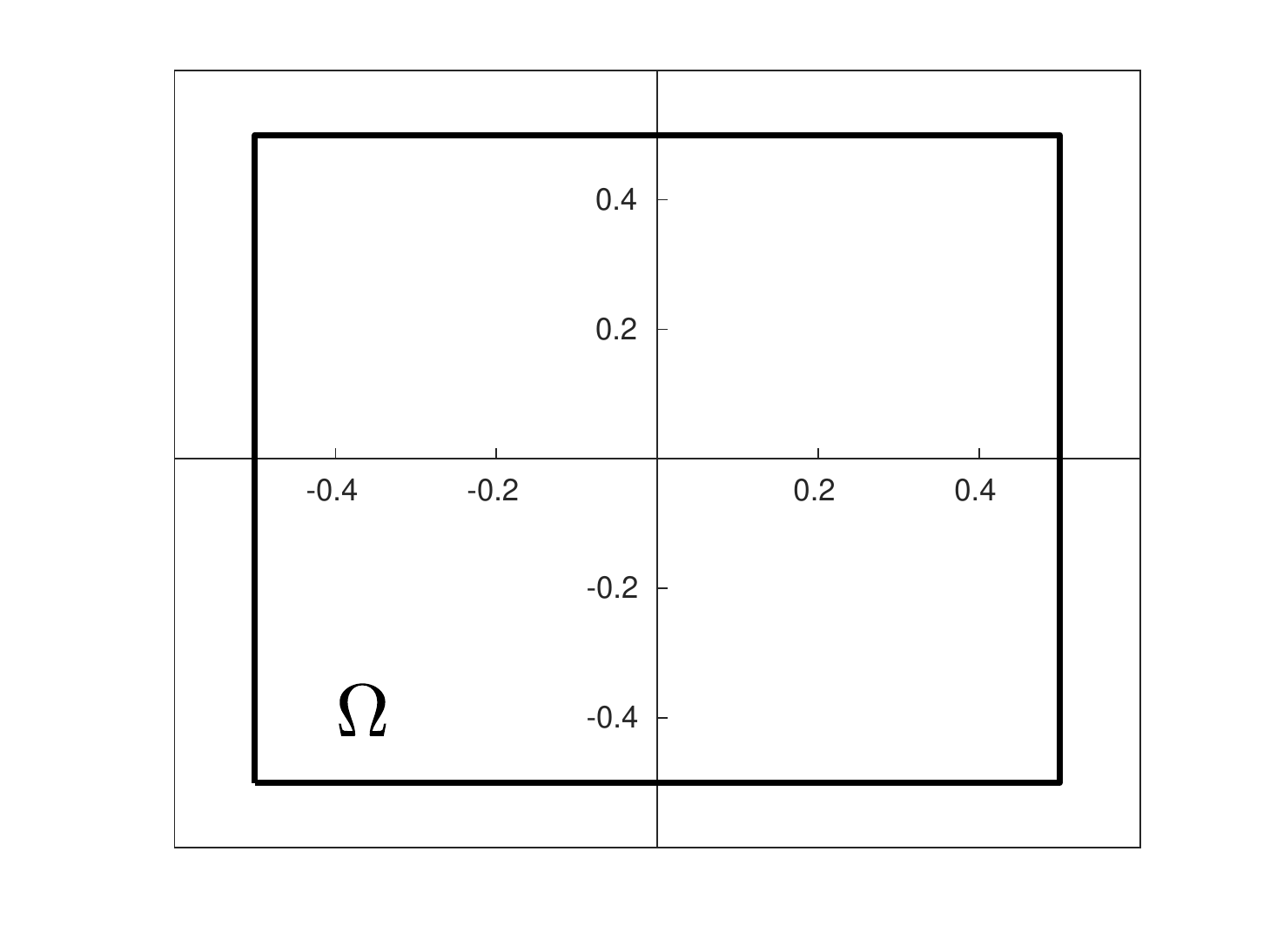}
\end{subfigure}
\hfill
\begin{subfigure}[t]{0.49\textwidth}
\includegraphics[width = \textwidth]{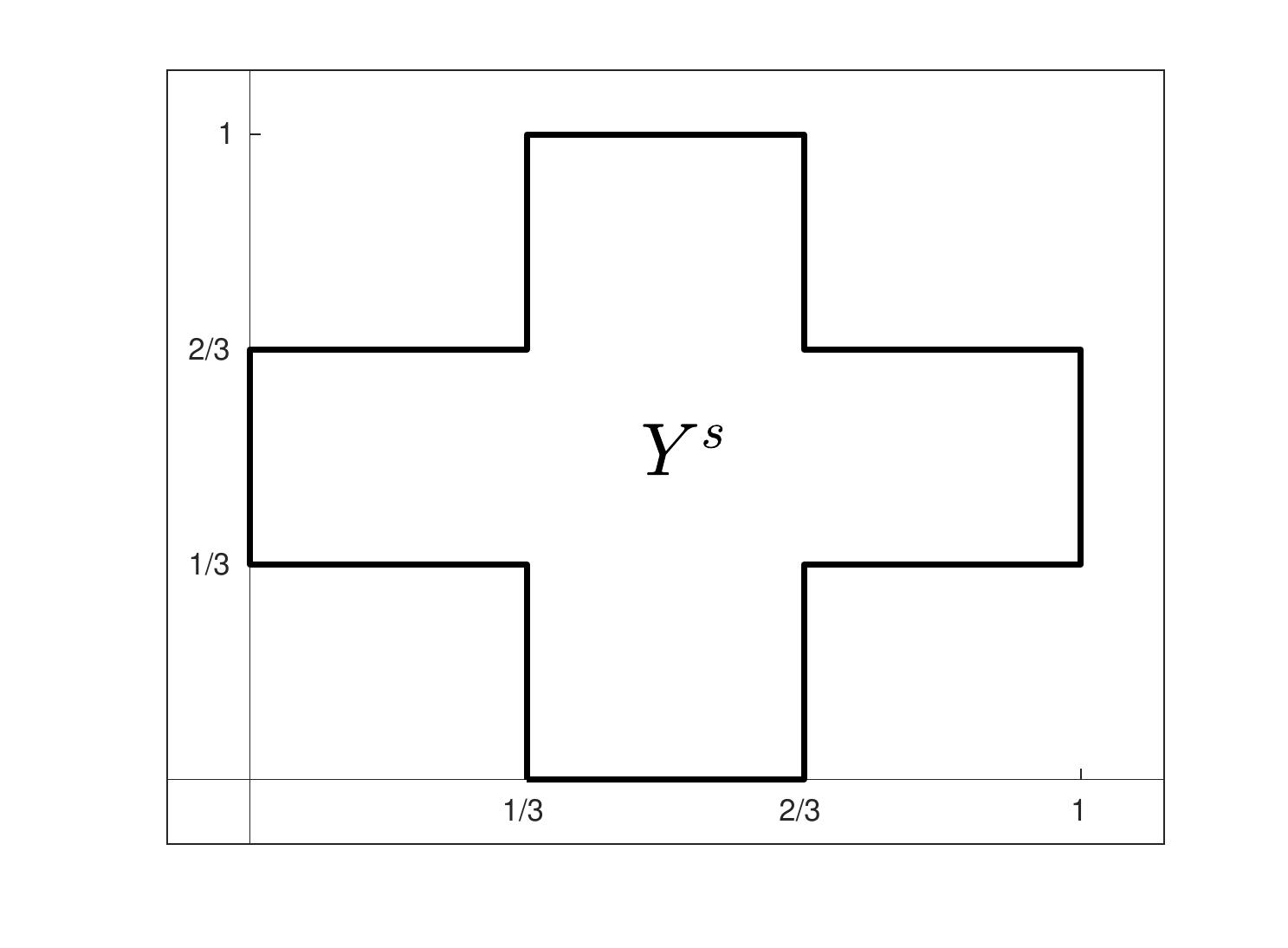}
\end{subfigure}
\caption{{\bf Non-dimensionalized domains} $\Omega = (-\frac{1}{2}, \frac{1}{2})^2$ for the effective, macroscopic elasticity-diffusion problem and $Y^s = \left((\frac{1}{3}, \frac{2}{3})\times (0,1) \right) \cup  \left( (0,1)\times (\frac{1}{3}, \frac{2}{3}) \right)$ for the elasticity and diffusion cell problems.}
\label{fig:domains}
\end{figure}
In the following, we focus on the simulation of the upscaled model summarized in section \ref{subsection:The_Upscaled_Equations_on_the_Fixed_Domain}. Due to the multi-scale character of the problem, computations have to be performed on two distinct domains: $\Omega$ for the macroscopic elasticity-diffusion system and $Y^s$ for the cell problems for elasticity and diffusion (see Figure \ref{fig:domains}). For the implementation, we have used the C$++$ finite element library \textit{deal.II} \cite{arndt2020deal} which provides comprehensive tools for mesh generation, dof handling, matrix assembly etc. All computations are performed on uniformly refined grids based on quadrilateral elements in two space dimensions. Spatial discretization for both, the macroscopic equations and the cell problems, is achieved using (bi-)linear \textit{Lagrange} finite elements. For temporal discretization of the effective, macroscopic diffusion equation $\eqref{homogenized_diffusion_problem:a}$--$\eqref{homogenized_diffusion_problem:d}$, the \textit{Crank-Nicolson} method is employed. Note that the effective, macroscopic elasticity problem $\eqref{eq:homogenized_elasticity_problem:a}$--$\eqref{eq:homogenized_elasticity_problem:c}$ is quasi-static, i.e. time dependence is only introduced due to the potentially time-dependent data $\mathbf{f}^\text{elast}$ and $\mathbf{h}$. Practically, this means we solve a series of static equations rather than a dynamic equation and therefore do not need to discretize a time derivative for this subproblem. We also emphasize that the system exhibits only a one-sided coupling. The effective coefficients for the diffusion equation, $\Jhom$ and $\Dhom$, depend on the solution of the elasticity problem in a nonlinear manner, but the elasticity problem can be solved independently from the diffusion problem. We exploit this structure to avoid having to solve a nonlinear system for the unknowns of the finite element scheme. Instead, in each time point $t^k$, we compute first the solution $\mathbf{u}^{k}$ to the elasticity problem, use this information to compute ${J^*}^k$ and ${\Dhom}^k$ and then solve the diffusion problem for $c^k$. For the treatment of the cell problems, one also has to differentiate: Since the elasticity cell problems do neither depend on the macroscopic equations nor on time or on the macroscopic space variable, it is sufficient to only solve them once for each $i,j = 1,2$ at the beginning of the simulation. The resulting effective elasticity tensor $\Ahom$ is constant in time and space during the subsequent simulation. In contrast, the diffusion cell problems depend on the solution of the macroscopic elasticity problems and therefore implicitly also on time and the macroscopic space variable through the coefficient $\mathbf{D}_0$ (see $\eqref{eq:J0_and_D0:b}$ and $\eqref{eq:diffusion_cell_problems:a}$). This means we have to solve different diffusion cell problems for each time point $t^k$ and each macroscopic spatial quadrature point $x_q$ of the chosen discretization, resulting in effective quantities $\Jhom$ and $\Dhom$ that change with time and space in the macroscopic diffusion equation. Let us underline again that this effect is introduced when the initial diffusion problem on the moving domain is transformed into an equation on a fixed domain in section \ref{subsection:Transformation_of_the_diffusion_equation_to_the_fixed_domain} and independent of the initial diffusion tensor $\widehat{\mathbf{D}}$ (cf. $\eqref{system:elasticity-diffusion-system_on_mixed_domains:e}$) being time- and/or space-dependent or not. The process of solving all the diffusion cell problems consumes most of the computation time during the simulation. Optimization of this step in terms of parallelization or some adaptive scheme is an interesting and promising task to save computation time but not within the scope of this paper. 
\subsection{The Model Problem}
\label{subsection:The_Model_Problem}
We fix now data for further investigations of the upscaled problem in section \ref{subsection:The_Upscaled_Equations_on_the_Fixed_Domain}. Let $\Omega = (-\frac{1}{2}, \frac{1}{2})^2$ be a quadratic domain and $Y^s = \left((\frac{1}{3}, \frac{2}{3})\times (0,1) \right) \cup  \left( (0,1)\times (\frac{1}{3}, \frac{2}{3}) \right)$ a cross-shaped domain (cf. Figure \ref{fig:domains}). For $\Gamma_D^\text{elast}$, we chose the lateral parts of $\partial \Omega$ and set $\Gamma_N^\text{elast} := \partial \Omega \setminus \Gamma_D^\text{elast}$. We do not consider a body force acting on the material and therefore set $\mathbf{f}^\text{elast} \equiv 0$. The deformation of the domain is solely induced by the time-periodic Dirichlet boundary condition
\begin{equation}
\mathbf{h}(t,x) = 
\begin{cases}
\left( a \frac{1 - \cos(2\pi f t)}{2}, 0 \right)^T & \text{ if } x  \in  \Gamma_D^\text{elast} \cap \{ x_1 = \frac{1}{2} \} \\ 
\left( -a \frac{1 - \cos(2\pi f t)}{2}, 0 \right)^T & \text{ if } x \in  \Gamma_D^\text{elast} \cap \{ x_1 = - \frac{1}{2} \},
\end{cases} \label{eq:elasticity_dirichlet_bc_constant_front}
\end{equation}
with $a = 0.25$, the maximal displacement at the boundary, and $f = 1$, the frequency. For an intuition about this boundary condition, see Figures \ref{fig:elasticity_visualisation:b} and \ref{fig:elasticity_visualisation:c}. We further assume that the material is isotropic, i.e. its elastic properties can be described by two \textit{Lam\'e}-constants $\lambda, \mu$ and the entries of the elasticity tensor are given by 
\begin{equation}
\mathbf{A}_{ijkl} = \mu ( \delta_{ik} \delta_{jl} + \delta_{il} \delta_{jk} ) +  \lambda \delta_{ij} \delta_{kl} \quad \text{ for } i,j,k,l = 1,2.
\end{equation}
For the non-dimensionalised elasticity tensor in the elasticity cell problems $\eqref{eq:elasticity_cell_problems:a}$, we set $\lambda, \mu = 1$. In Table \ref{tab:comparison_of_initial_and_effective_quantities}, we list the non-zero entries of $\mathbf{A}$ opposed to their counterparts in the effective elasticity tensor $\Ahom$ and see that the isotropy of $\mathbf{A}$ is not conserved during the upscaling process. For the diffusion problem, we set $\Gamma_D^\text{diff}$ to be the upper part of the boundary of $\Omega$ and $\Gamma_N^\text{diff} = \partial \Omega \setminus \Gamma_D^\text{diff}$. We do not consider a sink/source, i.e. we set $f^\text{diff} \equiv 0$ and prescribe a constant concentration $g \equiv 1$ at the Dirichlet boundary as well as no-flux conditions on the Neumann boundary together with initial condition $c^0 \equiv 0$. The non-dimensionalized and constant initial diffusion tensor in $\eqref{eq:J0_and_D0:b}$ is set to be
\begin{equation}
\widehat{\mathbf{D}} = 
\begin{pmatrix}
0.5 & 0 \\
0 & 0.5
\end{pmatrix}.
\label{eq:diffusion_coefficient_for_numerics}
\end{equation}
As mentioned before, the effective quantity $\Dhom$ is not constant in time or space in general, but depends on the displacement via the deformation gradient. Nonetheless, for $t = 0$, it is constant across the whole domain $\Omega$ since it is initially in the non-deformed state. Therefore, the deviation between $\widehat{\mathbf{D}}$ and $\Dhom(t = 0, \cdot)$ is only due to the heterogeneity of the domain and not due to deformation, when we compare the non-zero entries of these quantities in Table \ref{tab:comparison_of_initial_and_effective_quantities}. 
\begin{table}
\centering
\begin{tabular}{c c c c | c | c } 
 i & j & k & l & $\mathbf{A}$ & $\Ahom$\\
 \hline 
1 & 1 & 1 & 1 & 3 & 0.952656 \\
1 & 1 & 2 & 2 & 1 & 0.131924 \\
1 & 2 & 1 & 2 & 1 & 0.281493 \\
1 & 2 & 2 & 1 & 1 & 0.281493 \\
2 & 1 & 1 & 2 & 1 & 0.281493 \\
2 & 1 & 2 & 1 & 1 & 0.281493 \\
2 & 2 & 1 & 1 & 1 & 0.131924 \\
2 & 2 & 2 & 2 & 3 & 0.952656 \\
\hline
& & & & $\widehat{\mathbf{D}}$ & $\Dhom$ \\
\hline
1 & 1 & - & - & 0.5 & 0.184577 \\
2 & 2 & - & - & 0.5 & 0.184577
\end{tabular}
\caption{{\bf Comparison} of the non-zero entries of the initial elasticity and diffusion tensors $\mathbf{A}$ and $\widehat{\mathbf{D}}$ with their effective counterparts $\Ahom$ and $\Dhom$. As $\Dhom$ is in general a function of time and space, we display it here only for $t = 0$. $\Dhom$ is then constant in the $x$-variable since the domain is initially not deformed at all.}
\label{tab:comparison_of_initial_and_effective_quantities}
\end{table}
\subsection{Numerical Convergence Studies}
\label{subsection:Numerical_Convergence_Studies}
\begin{table}
\centering
\resizebox{\textwidth}{!}{\begin{tabular}{ c r c | c c c c | c c c c }
 \multicolumn{1}{c}{} &\multicolumn{1}{c}{} &\multicolumn{1}{c}{} & \multicolumn{4}{c}{\textbf{concentration}} & \multicolumn{4}{c}{\textbf{displacement}} \\
cycle & \#cells & $h$ & $L^2$ & EOC($L^2$) & $H^1$ & EOC($H^1$) & $L^2$ & EOC($L^2$) & $H^1$ & EOC($H^1$)  \\ 
\hline \hline
0 & 1 & 1.414e$+$00 & - & - & - & - & - & - & - & - \\
1 & 4 & 7.071e$-$01 & 1.146e$-$01 & - & 4.375e$-$01 & - & 1.982e$-$02   & - & 8.172e$-$02    & - \\
2 & 16 & 3.536e$-$01 & 3.819e$-$02 & 1.59 & 2.466e$-$01 & 0.83 &  4.707e$-$03 & 2.07 & 3.832e$-$02 & 1.09 \\
3 & 64 & 1.768e$-$01 & 1.039e$-$02 & 1.88 & 1.470e$-$01 & 0.75 &  1.518e$-$03 & 1.63 & 2.317e$-$02 & 0.73 \\
4 & 256 & 8.839e$-$02 & 3.639e$-$03 & 1.51 & 9.388e$-$02 & 0.65 &  5.113e$-$04 & 1.57 & 1.440e$-$02 & 0.69 \\
5 & 1024 & 4.419e$-$02 & 1.341e$-$03 & 1.44 & 6.576e$-$02 & 0.51 & 1.746e$-$04& 1.55 & 8.892e$-$03 & 0.70 \\
6 & 4096 & 2.209e$-$02 & 3.887e$-$04 & 1.79 & 3.747e$-$02 & 0.81 & 5.854e$-$05 & 1.58 & 5.389e$-$03 & 0.72 \\
7 & 16384 & 1.105e$-$02 &  1.033e$-$04 & 1.91 & 1.950e$-$02 & 0.94 &  1.914e$-$05 & 1.61 & 3.209e$-$03 & 0.75 
\end{tabular}}
\resizebox{\textwidth}{!}{\begin{tabular}{c r c | c c c c | c c c c}
  \multicolumn{1}{c}{} &\multicolumn{1}{c}{} &\multicolumn{1}{c}{}  & \multicolumn{4}{c}{\textbf{concentration}} & \multicolumn{4}{c}{\textbf{displacement}} \\
cycle & \#cells & $h$ & $L^2$ & EOC($L^2$) & $H^1$ & EOC($H^1$) & $L^2$ & EOC($L^2$) & $H^1$ & EOC($H^1$)  \\
\hline \hline
0 & 1 & 1.414e$+$00 & - & - & - & - & - & - & - & -  \\
1 & 4 & 7.071e$-$01 & 1.104e$-$01 & - & 4.554e$-$01 & - &  1.250e$-$01 & - &  5.153e$-$01 & - \\
2 & 16 & 3.536e$-$01 & 5.234e$-$02 & 1.08 & 2.947e$-$01 & 0.63 & 2.630e$-$02 & 2.25 & 1.872e$-$01 & 1.46 \\
3 & 64 & 1.768e$-$01 & 1.381e$-$02 & 1.92 & 1.598e$-$01 & 0.88 & 6.306e$-$03 & 2.06 & 8.419e$-$02 & 1.15 \\
4 & 256 & 8.839e$-$02 & 3.861e$-$03 & 1.84 & 7.927e$-$02 & 1.01 & 1.580e$-$03 & 2.00 & 4.111e$-$02 & 1.03 \\
5 & 1024 & 4.419e$-$02 & 1.001e$-$03 & 1.95 & 3.996e$-$02 & 0.99 &  3.979e$-$04 & 1.99 & 2.049e$-$02 & 1.00 \\
6 & 4096 & 2.209e$-$02 & 2.545e$-$04 & 1.98 & 2.010e$-$02 & 0.99 & 1.000e$-$04 & 1.99 & 1.025e$-$02 & 1.00 \\
7 & 16384 & 1.105e$-$02 & 6.419e$-$05 & 1.99 & 1.009e$-$02 & 0.99 &  2.508e$-$05 & 2.00 & 5.132e$-$03 & 1.00
\end{tabular}}
\caption{{\bf Errors and EOC for the displacement and the concentration} with respect to the discretization parameter $h$ of the grid for the macroscopic domain $\Omega$. The computations have been performed for the model problem with mixed boundary conditions as described in Section \ref{subsection:The_Model_Problem} (upper table) and for the model problem with pure Dirichlet boundary conditions as described in Section \ref{subsection:Numerical_Convergence_Studies} (lower table).}
\label{tab:convergence}
\end{table}
\begin{figure}
\centering
\begin{subfigure}[t]{0.49\textwidth}
\includegraphics[width = \textwidth]{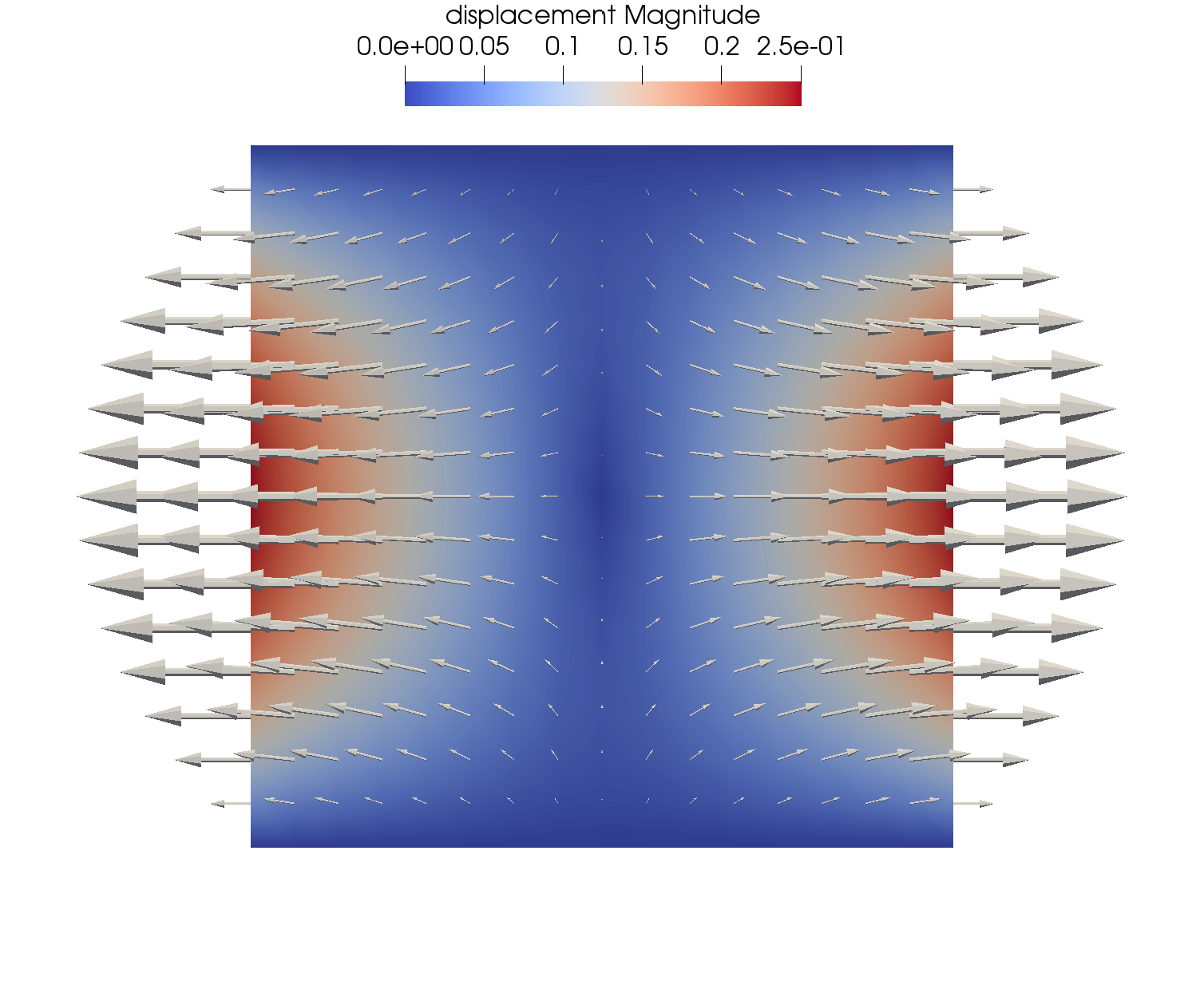}
\end{subfigure}
\hfill
\begin{subfigure}[t]{0.49\textwidth}
\includegraphics[width = \textwidth]{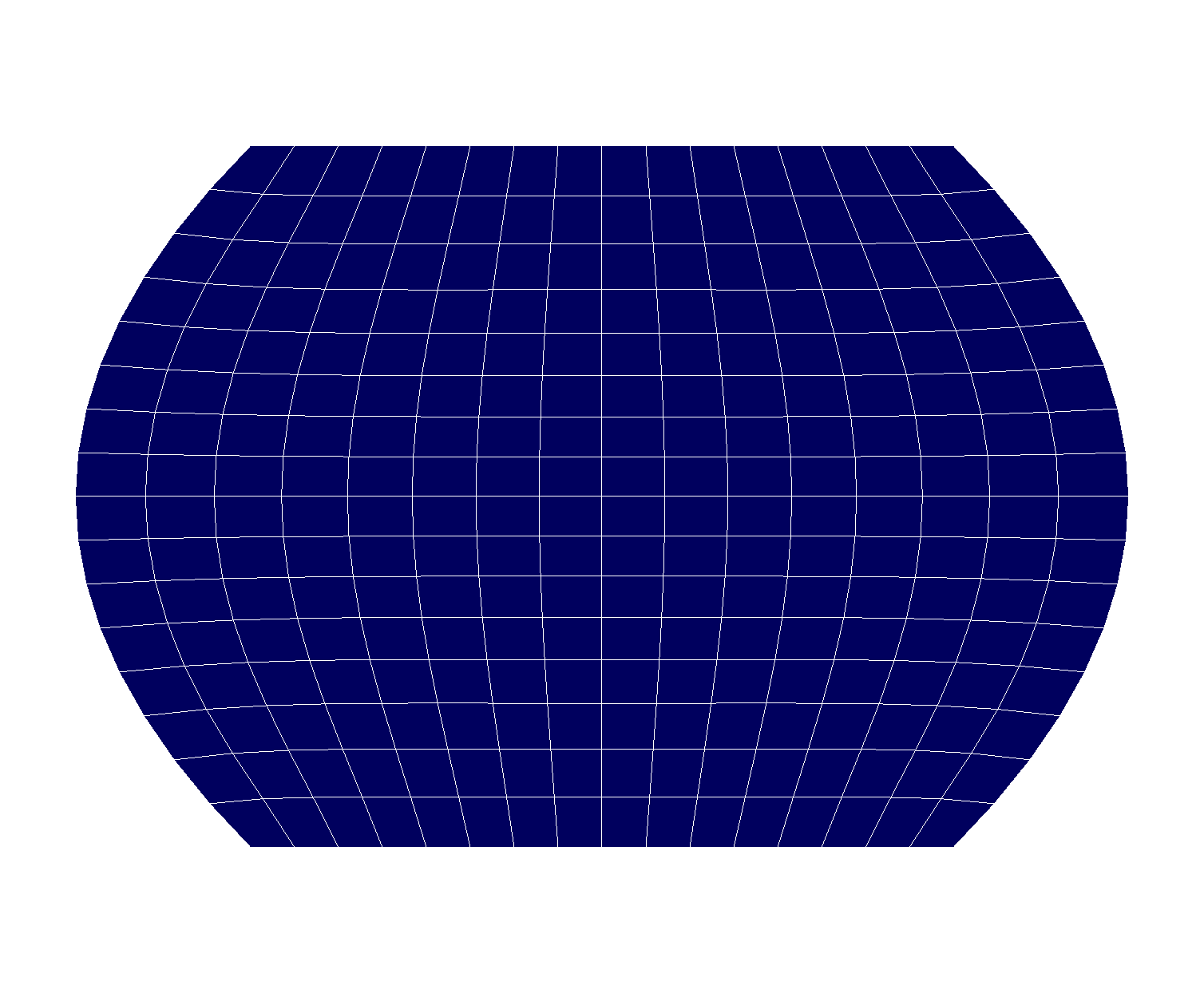}
\end{subfigure}
\caption{{\bf The displacement for the modified model problem} with pure Dirichlet boundary condition. The problem is used for the investigations concerning the convergence properties of solutions.}
\label{fig:parabola_plots}
\end{figure}
\begin{figure}
\centering
\begin{subfigure}[t]{0.49\textwidth}
\includegraphics[width = \textwidth]{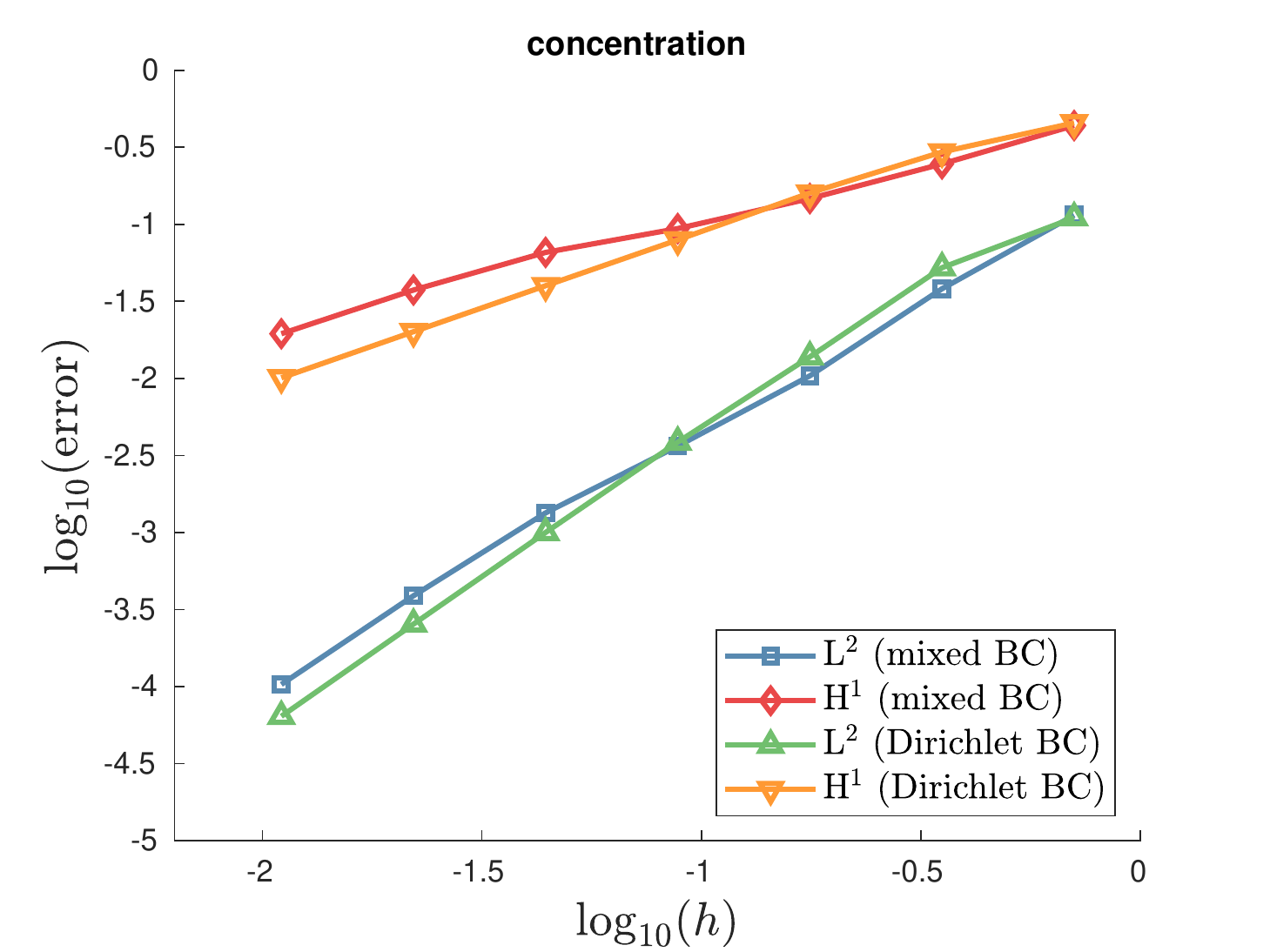}
\end{subfigure}
\hfill
\begin{subfigure}[t]{0.49\textwidth}
\includegraphics[width = \textwidth]{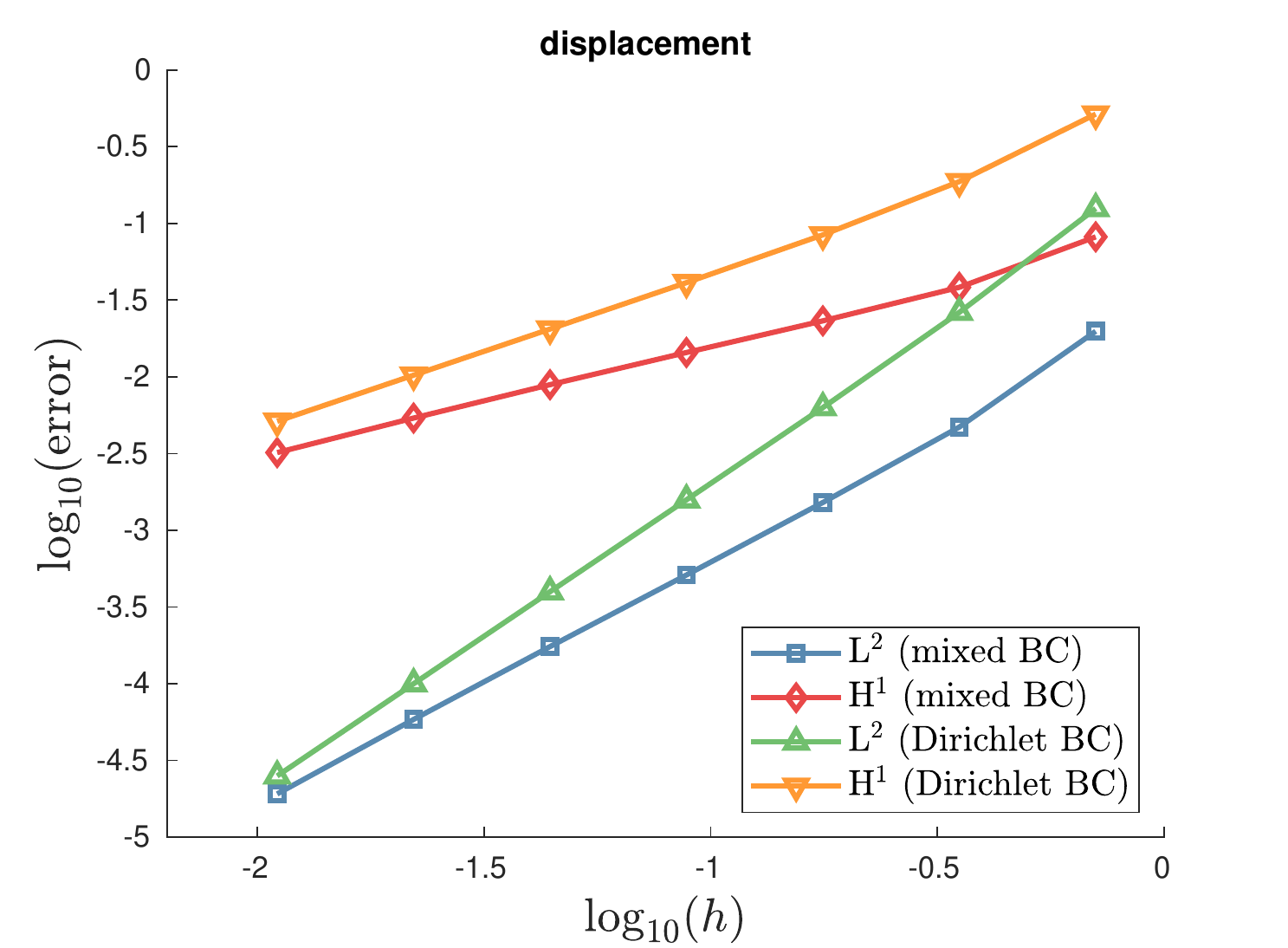}
\end{subfigure}
\caption{{\bf Convergence plots} for the data given in Table \ref{tab:convergence}. We compare convergence properties of the model problem with mixed boundary conditions described in Section \ref{subsection:The_Model_Problem} and the model problem with pure Dirichlet boundary condition described in Section \ref{subsection:Numerical_Convergence_Studies} with respect to the discretization parameter $h$ of the grid for the macroscopic domain $\Omega$.}
\label{fig:convergence_plots}
\end{figure}
Before we proceed with a more detailed view on the solutions of the proposed model problem, let us first investigate the convergence of the finite element scheme. Due to the two-scale character of the problem, there are two distinct computational domains, representing $\Omega$ and $Y^s$, respectively. In the following, we restrict ourselves to the analysis of the macroscopic solutions $\mathbf{u}$ and $c$ with respect to the discretization parameter $h$ of the grid for $\Omega$, where $h$ is the diameter of the largest quadrilateral in the triangulation.  For the simulations, we use the model problem presented in the previous section \ref{subsection:The_Model_Problem}. To quantify its convergence properties, we compute the error between solutions at the fixed time $t = 1.5$ on subsequently refined grids, using an uniform refinement strategy to obtain finer meshes. Additionally, we compute the estimated order of convergence (EOC) according to the following formulas:
\begin{equation}
\text{error}_i = \| \square^{i-1} - \square^{i} \|_*, \quad \text{EOC}_i = \log_{2} \left( \frac{\text{error}_{i-1}}{\text{error}_i} \right),
\label{eq:error_and_EOC}
\end{equation}
where $i = 1,2,...$ indicates the refinement cycle, $\square \in \{ \mathbf{u}, c \}(t = 1.5, \cdot)$, and $* \in \{L^2(\Omega), H^1(\Omega) \}$. The so-obtained data is listed in the upper part of Table \ref{tab:convergence}. Note that the EOCs for the most part lie below the convergence orders known from theory for the implemented finite element method of linear lagrangian elements (order two for $L^2$-norm and order one for $H^1$-norm). This suboptimal behavior might be attributed to the mixed boundary conditions in our model problem, which are known to be a potential source of decreased regularity of solutions. In fact, we are able to recover the expected convergence rates, if we simulate instead a related problem with pure Dirichlet boundary condition and matching initial condition for the diffusion problem. Let $\Gamma_D^\text{elast} = \Gamma_D^\text{diff} = \partial \Omega$, $\Gamma_D^\text{elast} = \Gamma_D^\text{diff} = \emptyset$ and set $g \equiv 1$ on $\Gamma_D^\text{diff}$, $c^0 \equiv 1$ and
\begin{equation*}
\mathbf{h}(t,x) = 
\begin{cases}
\left( \pm a\frac{1 - \cos(2 \pi f t)}{2} \frac{0.25 - x_2^2}{0.25}, 0 \right)^T, & x \in \Gamma_D^\text{elast} \cap \{ x_1 = \pm \frac{1}{2} \}  \\
(0,0)^T,  & x \in \Gamma_D^\text{elast} \cap \{ x_2 = \pm \frac{1}{2} \} ,
\end{cases}
\end{equation*}
with, as above, $a = 0.25$, the maximal displacement at the boundary and $f = 1$, the frequency. This boundary condition for the elasticity problem results in a parabola-shaped extension of the domain at the lateral parts of the boundary while the upper and lower parts of the boundary is clamped, see Figure \ref{fig:parabola_plots}. Analogously to the study of the problem with mixed boundary conditions, we gather convergence data for simulations of this problem in the lower part of Table \ref{tab:convergence}. Comparison clearly indicates that the loss of convergence speed comes from the imposition of the mixed boundary conditions. However, the actual error between subsequent solutions appears to be smaller in the case of the problem with mixed boundary condition, at least for courser meshes, see Figure \ref{fig:convergence_plots}. Of course, the opposite is eventually true for finer meshes.
\subsection{The Solution of the Elasticity Subproblem}
\begin{figure}[t]
\centering
\begin{subfigure}[b]{0.32\textwidth}
\includegraphics[width = \textwidth]{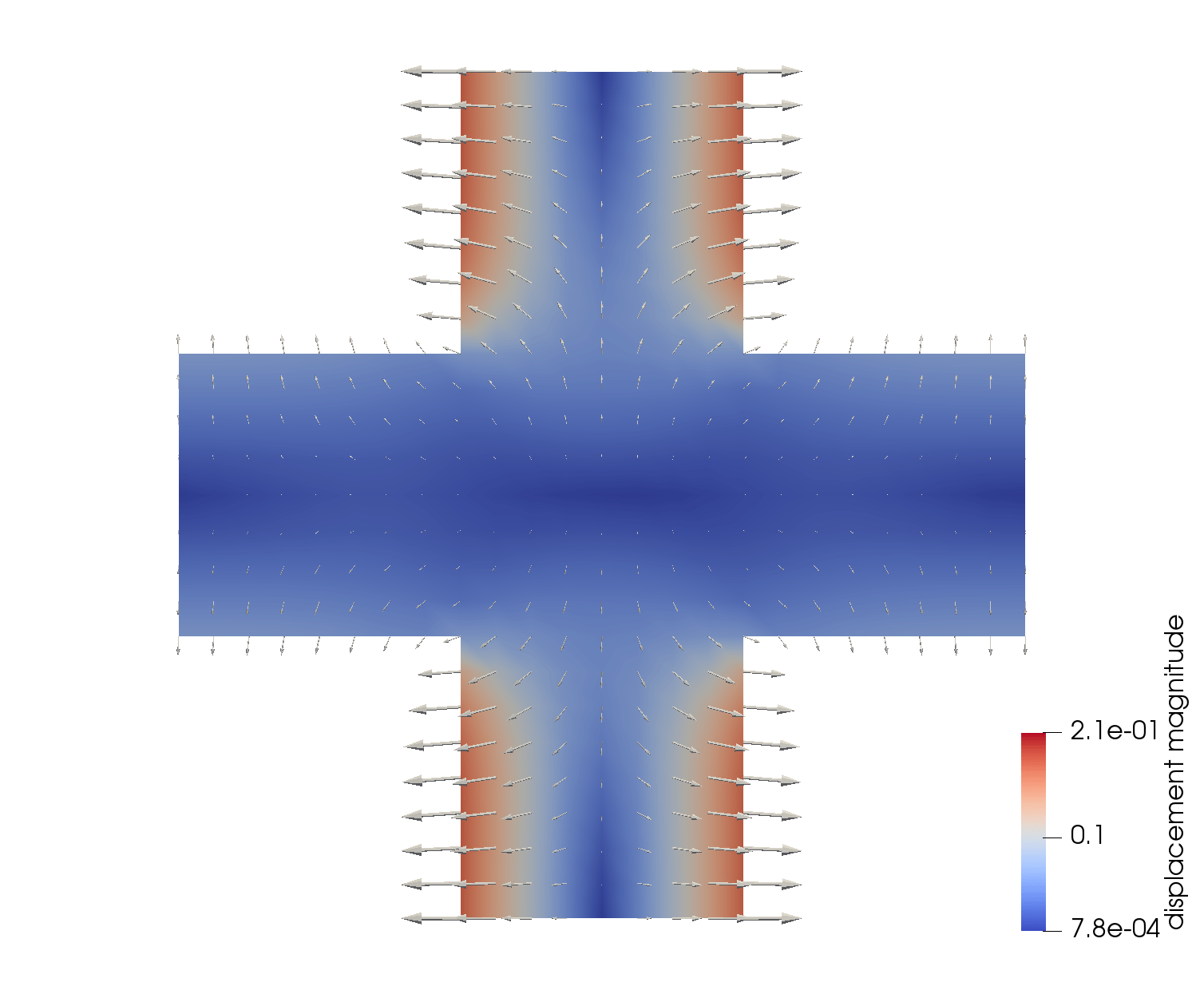}
\caption{$\boldsymbol{\chi}_{11}$}
\end{subfigure}%
\begin{subfigure}[b]{0.32\textwidth}
\includegraphics[width = \textwidth]{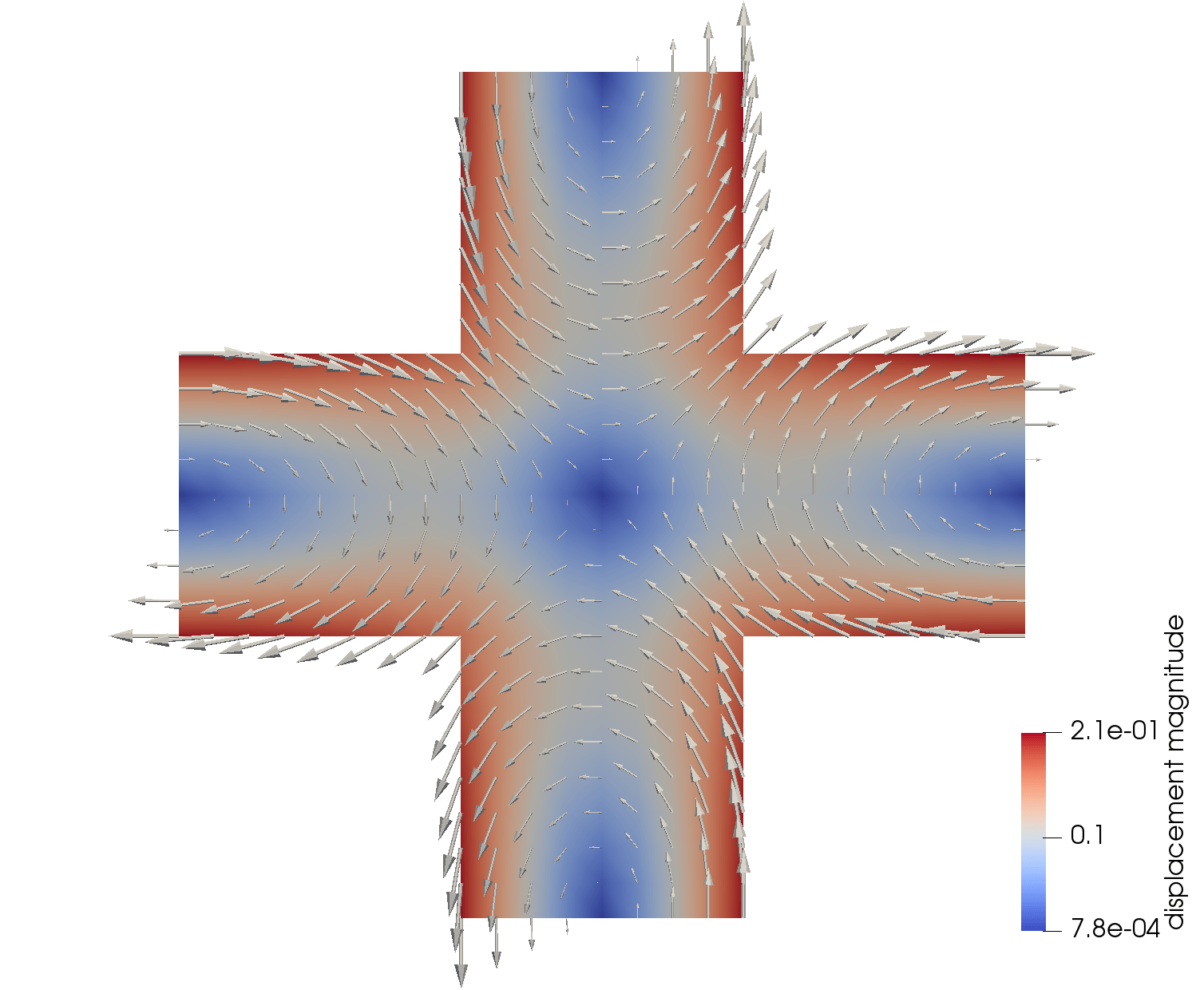}
\caption{$\boldsymbol{\chi}_{12} = \boldsymbol{\chi}_{21}$}
\end{subfigure}%
\begin{subfigure}[b]{0.32\textwidth}
\includegraphics[width = \textwidth]{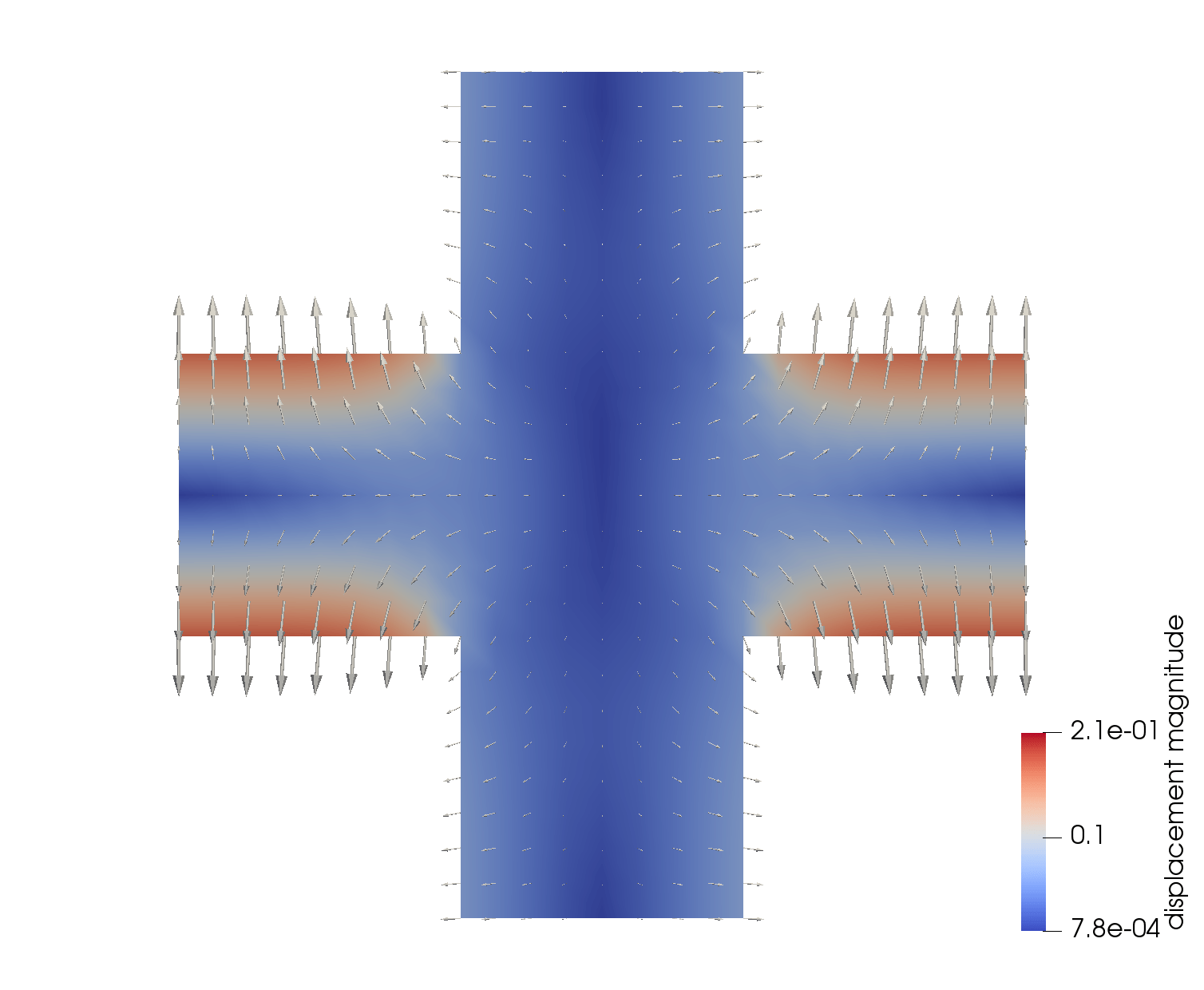}
\caption{$\boldsymbol{\chi}_{22}$}
\end{subfigure}
\caption{{\bf Cell solutions} of the elasticity cell problems $\eqref{eq:elasticity_cell_problems}$. Note that the symmetry of $\mathbf{A}$ implies $\boldsymbol{\chi}_{12} = \boldsymbol{\chi}_{21}$.}
\label{fig:elasticity_cell_solutions}
\end{figure}
\begin{figure}[t]
\begin{minipage}[b]{.65\textwidth}
  \centering
  \includegraphics[height = 5cm]{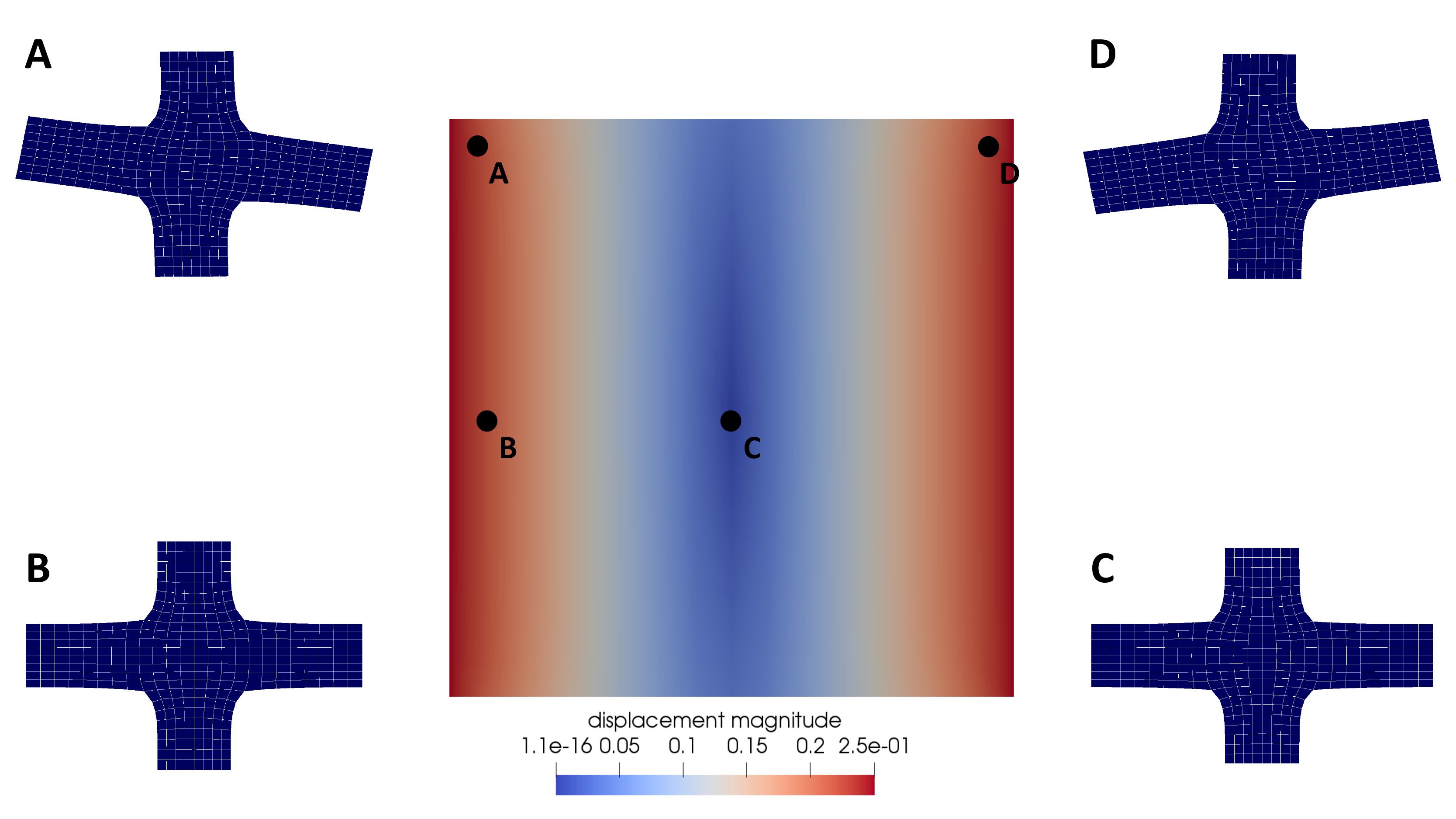}
  \subcaption{Visualization of the magnitude of displacement $\mathbf{u}(t = 0.5, \cdot)$ on $\Omega$ together with deformed microscopic domains $Y^s$ for different points in $\Omega$.}
  \label{fig:elasticity_visualisation:a}
\end{minipage}%
\hfill
\begin{minipage}[b]{.33\textwidth}
  \centering
  \includegraphics[height = 2.2cm]{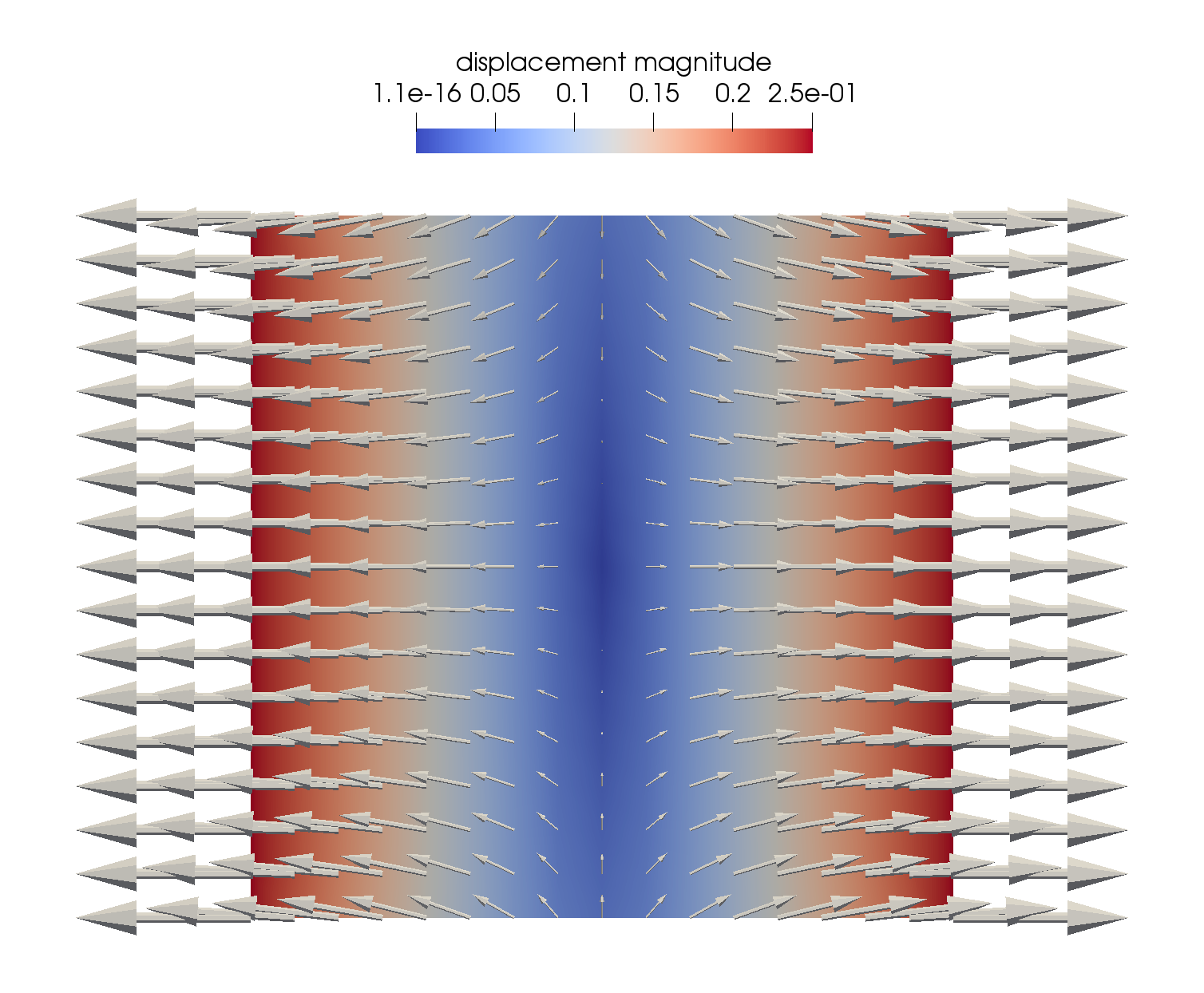}
  \subcaption{The displacement $\mathbf{u}(t=0.5, \cdot)$.}
  \label{fig:elasticity_visualisation:b}
  \includegraphics[height = 2.2cm]{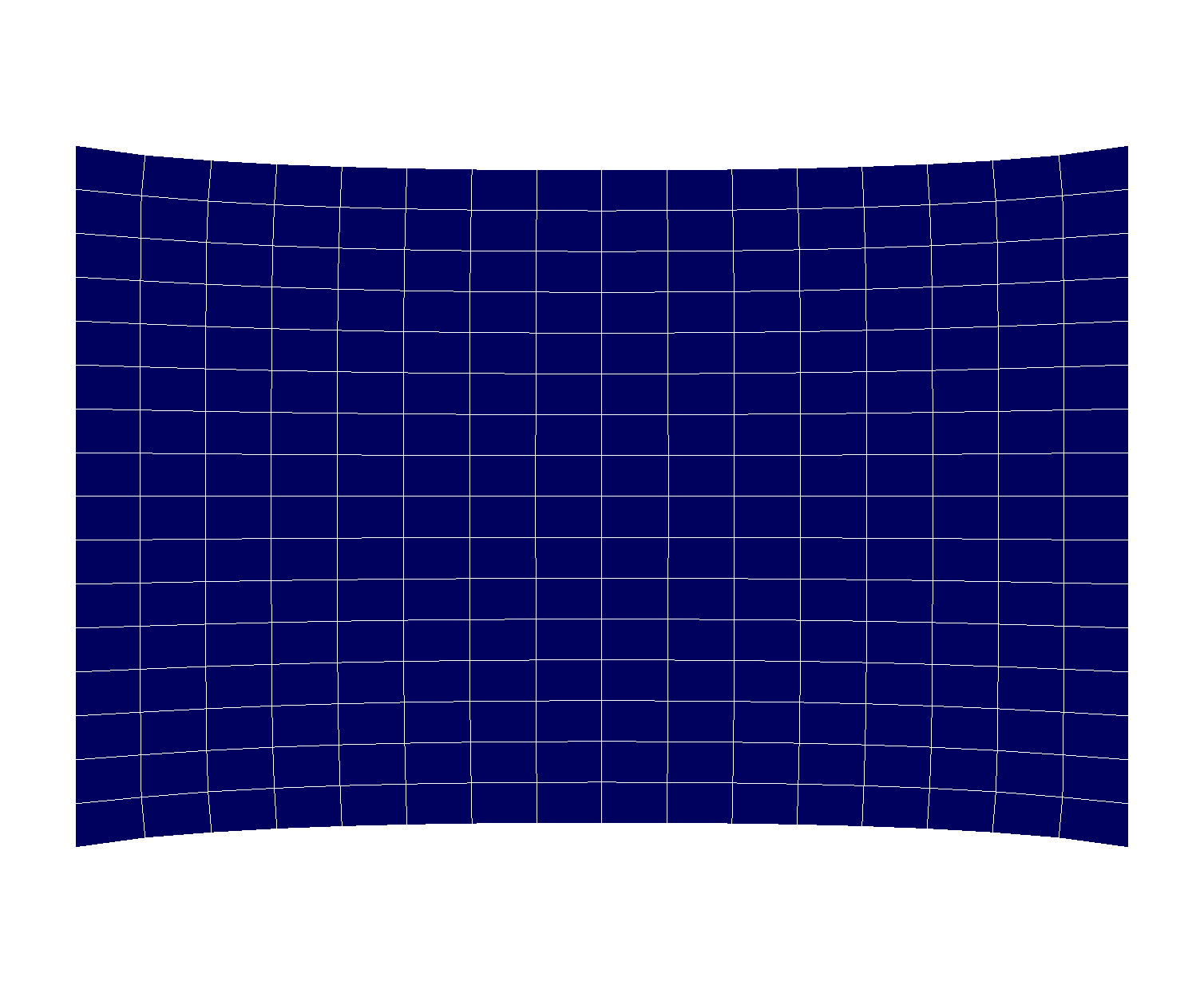}
  \subcaption{The shape of $\Omega^0(t = 0.5)$.}
  \label{fig:elasticity_visualisation:c}
\end{minipage}
\caption{{\bf Visualizations associated to the solution of the elasticity subproblem.} In a) and b), the magnitude of the displacement is encoded as color of the domain. In b) the arrow size is proportional to the magnitude of the displacement.}
\label{fig:elasticity_visualisation}
\end{figure}
Let us now give a more detailed analysis of the numerical solution to the effective elasticity subproblem using the aforementioned data. As stated earlier, it suffices to solve the elasticity cell problems $\eqref{eq:elasticity_cell_problems}$ once at the beginning of the simulation due to the circumstance of the coefficient $\mathbf{A}$ being constant. Note that the symmetry properties of $\mathbf{A}$ imply the symmetry $\boldsymbol{\chi}_{12} = \boldsymbol{\chi}_{21}$. The resulting cell solutions $\boldsymbol{\chi}_{ij}$, $i,j = 1,2$, contain the information needed for the computation of the effective elasticity tensor $\Ahom$ in $\eqref{Ahom}$. For a visualization of the elasticity cell solution, see Figure \ref{fig:elasticity_cell_solutions}. The non-zero entries of $\Ahom$ are given in Table \ref{tab:comparison_of_initial_and_effective_quantities}. The cell solutions are also contained in the representation of the first order term $\mathbf{u}_1$, see $\eqref{eq:u1_representation}$, in the expansion $\eqref{eq:two-scale-expansion}$ for $\ue$. In fact, we are now equipped to give an approximation of $\ue$ in terms of the asymptotic expansion up to terms of order $O(\eps^2)$ and higher. To this end, let us investigate the $m$-th component of $\ue = (u_\eps^1, u_\eps^2)^T$. For a fixed point $x \in \Oes$, we write $x = x_* + \eps y$ with $x_* = \eps \mathbf{k}$, $\mathbf{k} \in \mathbb{Z}^2$, $y \in Y^s$, and obtain by using Taylor Expansion and the symmetry $\boldsymbol{\chi}_{12} = \boldsymbol{\chi}_{21}$
\begin{align*}
u_\eps^m (t,x) &= u^m(t, x_* + \eps y) + \eps \sum_{i,j = 1}^2 e(\mathbf{u})_{ij}(t, x_* + \eps y) \chi^m_{ij}(y) + O(\eps^2) \\
&= u^m(t, x_*) + \eps \left( \sum_{j = 1}^2 \frac{\partial u^m(t,x_*)}{\partial x_j}y_j + \sum_{i,j = 1}^2 \frac{\partial u^i(t, x_*)}{\partial x_j}\chi^m_{ij}(y)	\right) + O(\eps^2) \\
&= u^m(t, x_*) + \eps \left( \sum_{i,j = 1}^2 \frac{\partial u^i(t,x_*)}{\partial x_j} \left( y_j \delta_{im} + \chi^m_{ij}(y) \right) 	\right) + O(\eps^2).
\end{align*}
In Figure \ref{fig:elasticity_visualisation:a}, we give a visualization of the $\eps$-term in the above computation for $y \in Y^s$ at $t = 0.5$, i.e. when the material is maximally extended, to illustrate how the deformation acts on individual, microscopic crosses $Y^s$ which constitute $\Oes$. Additionally, the effective displacement $\mathbf{u}(t=0.5, \cdot)$ is displayed in Figure \ref{fig:elasticity_visualisation:b}. There, the displacement $\mathbf{u}$, defined on $\Omega$, is represented by arrows with size proportional to their magnitude. We can also use $\mathbf{u}$ to visualize the deformed domain $\Omega^0(t) := \{ \widehat{x} \in \mathbb{R}^2  \mid  \widehat{x} = x + \mathbf{u}(t,x), \; x \in \Omega \}$ as can be seen in Figure \ref{fig:elasticity_visualisation:c} for time $t = 0.5$. This gives us also a clear intuition for the effect of the boundary conditions of the elasticity problem stated earlier: Within one unit of time, the domain is stretched in lateral direction to attain the shape depicted in \ref{fig:elasticity_visualisation:c}. Then, the displacement at the boundary is released until the domain is again in the non-deformed configuration.
\subsection{The Solution of the Diffusion Subproblem}%
\label{subsection:Solution_of_the_Diffusion_Subproblem}
\begin{figure}[t]
\centering
\begin{subfigure}[b]{0.49\textwidth}
\includegraphics[width = \textwidth]{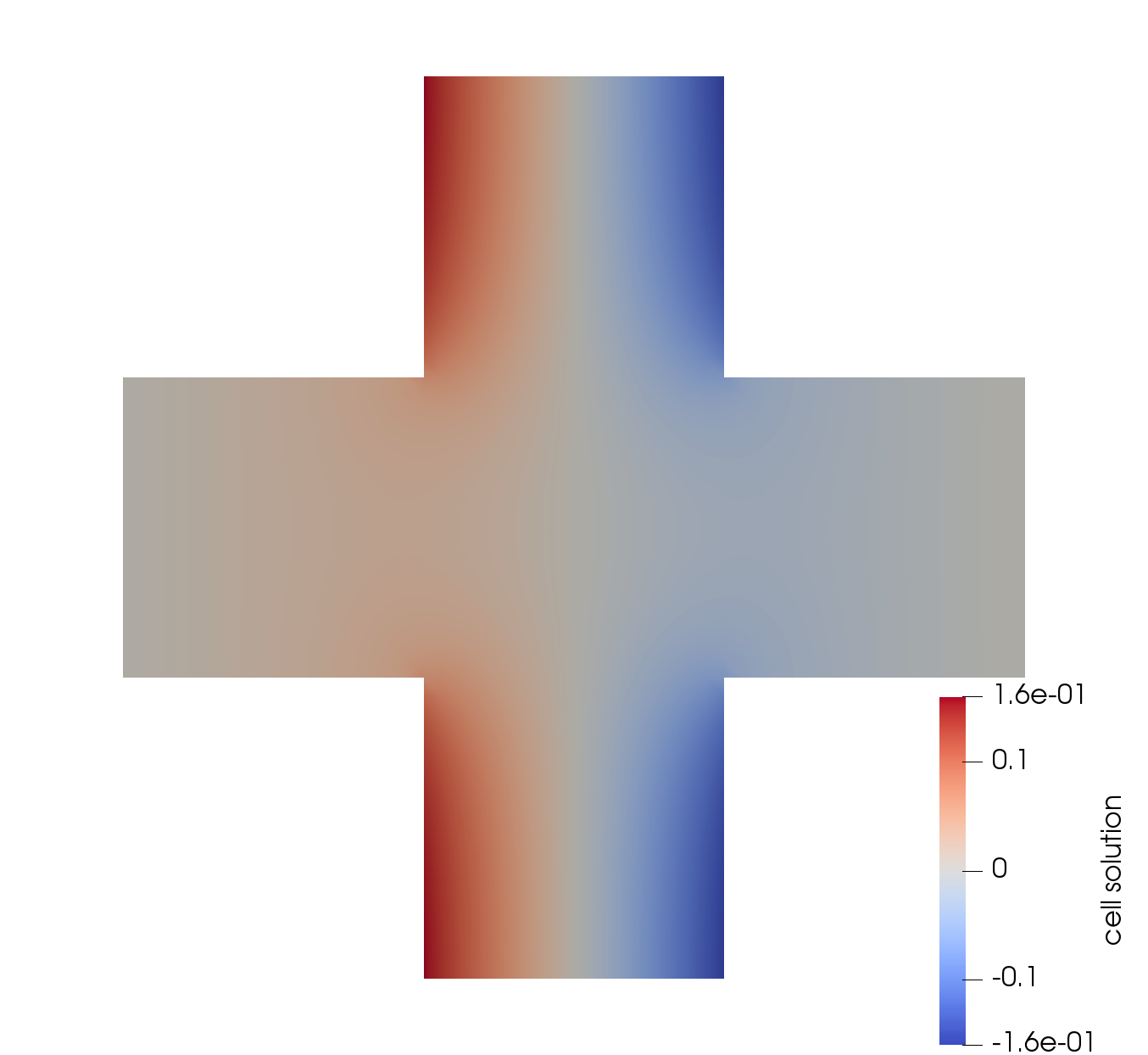}
\caption{$\eta_1(t = 0, x = (0,0)^T, y)$}
\end{subfigure}%
\begin{subfigure}[b]{0.49\textwidth}
\includegraphics[width = \textwidth]{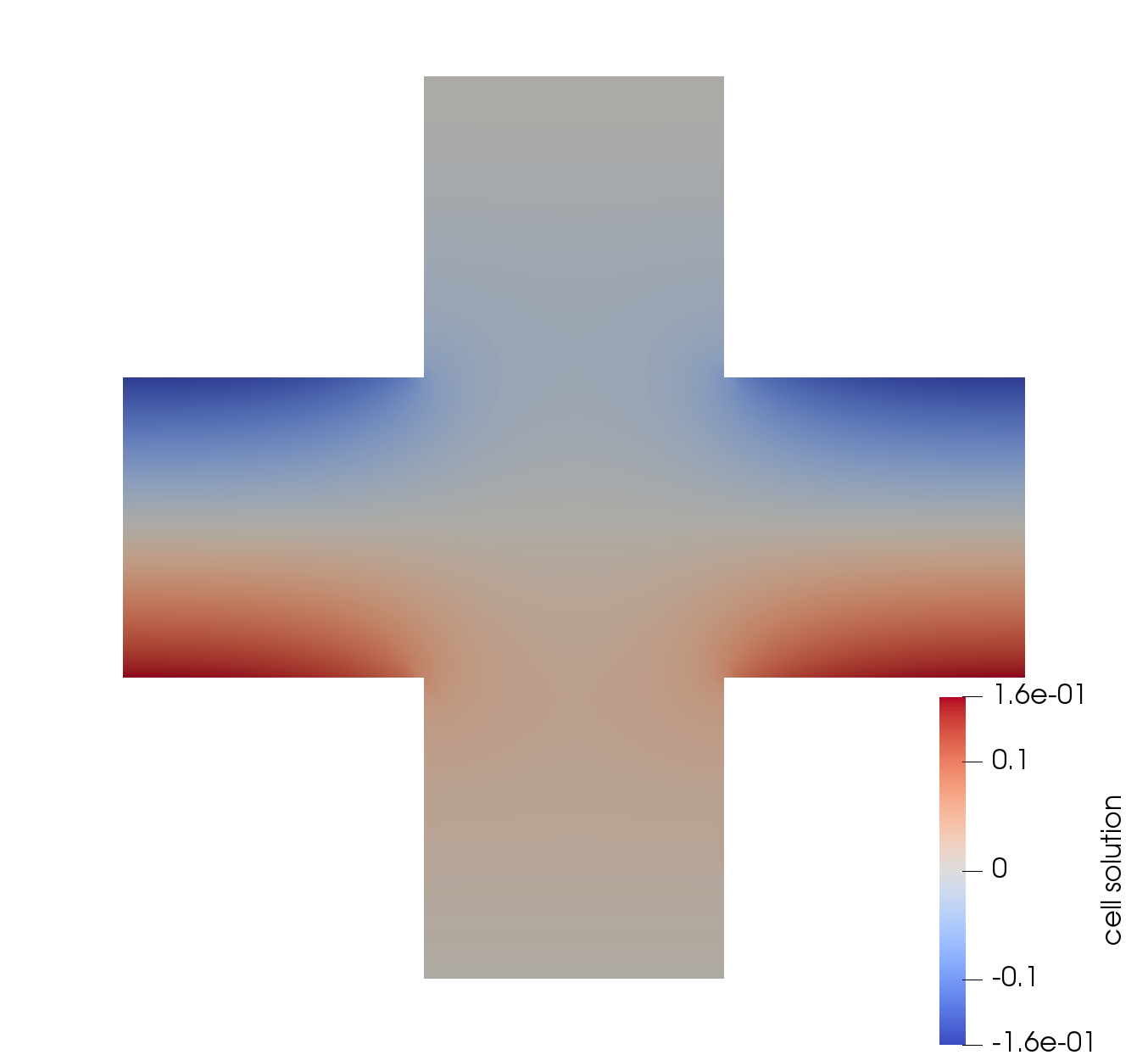}
\caption{$\eta_2(t = 0, x = (0,0)^T, y)$}
\end{subfigure}%
\caption{{\bf Cell solutions} of the diffusion cell problems $\eqref{eq:diffusion_cell_problems}$ for fixed $t = 0$ and $x = (0,0)^T$.}
\label{fig:diffusion_cell_solutions}
\end{figure}
\begin{figure}[t]
\centering
\includegraphics[width = 0.75\textwidth]{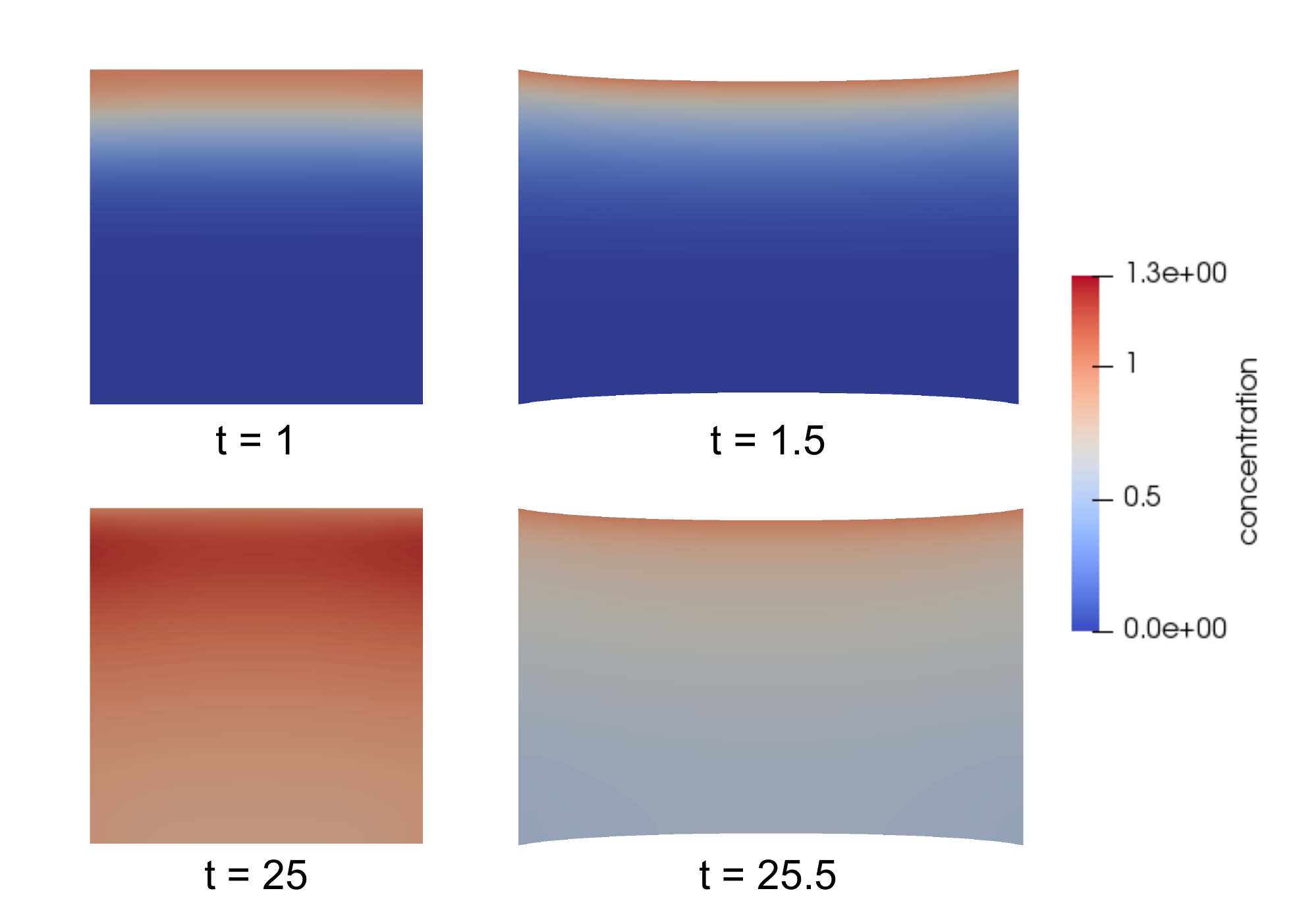}
\caption{{\bf The macroscopic concentration} $\widehat{c}$, visualized on the current deformed domain $\Omega^0(t)$ during simulations of the model problem from Section \ref{subsection:The_Model_Problem}.}
\label{fig:Macroscopic_Diffusion}
\end{figure}
Finally, we investigate the solution of the diffusion equation for the model problem from Section \ref{subsection:The_Model_Problem}. As mentioned before, the diffusion problem depends on the elasticity problem through its coefficients, but not vice versa. Therefore, in each time step, we solve first for the displacement, compute from this the effective coefficients $\Dhom$, $\Jhom$ for the diffusion problem and then solve the diffusion problem. The effective coefficients for the diffusion problem are defined by means of the diffusion cell solutions $\eta_i$ and the coefficient $\mathbf{D}_0$ for the diffusion cell problems is given in terms of $\mathbf{F}_0$, i.e. the lowest order term in the expansion of the deformation gradient $\mathbf{F}_\eps$. Concerning the numerical scheme, this amounts to \textit{dim} $\times$ \textit{number of time steps} $\times$ \textit{number of elements} $\times$ \textit{number of quadrature points per element} diffusion cell problems that have to be solved numerically during the simulation. For illustration, numerical approximations of the functions $\eta_i(t = 0.5, x = (0,0)^T, y)$, $i = 1,2$, $y \in Y^s$, are plotted in Figure \ref{fig:diffusion_cell_solutions}. Unsurprisingly, the computation of the diffusion cell solutions takes up most of the computation time, but there is some potential for optimization: The diffusion cell problems are independent from each other, so their treatment can be parallelized. Additionally, adaptive schemes such as clustering "similar" cell problems and solving only one representative cell problem per cluster are conceivable and have been successfully applied in the context of multi-scale models, see e.g. \cite{bastidas2020adaptive, garttner2020efficiency}.\\ Due to the transformation of the diffusion problem onto a fixed domain in Section \ref{subsection:Transformation_of_the_diffusion_equation_to_the_fixed_domain}, the concentration $c$, defined on $\Omega$, does not have a direct physical interpretation, as concentrations in the real world would be associated with points in the current deformed domain. Nonetheless, the concentration $\widehat{c}$, defined on the current deformed domain, can be obtained from the concentration $c$ using the following transformation
\begin{gather*}
\widehat{c}(t,\hx) := c(t, \mathbf{S}_0^{-1}(t, \hx)) \; \text{for} \; \hx \in \Omega^0(t) \; \text{with} \; \mathbf{S}_0(t,x) := x + \mathbf{u}\tx, \\ \text{and} \; \Omega^0(t) := \{ \hx \in \mathbb{R}^n, \; \widehat{x} = \mathbf{S}_0\tx, \; x \in \Omega \}.
\end{gather*}
The concentration $\widehat{c}$ from the simulations to the model problem from Section \ref{subsection:The_Model_Problem} is depicted for different time points in Figure \ref{fig:Macroscopic_Diffusion}. There, one can also spot an interesting effect of the time-periodic deformation of the domain on the accumulation of particles in the domain: When the domain is expanded, the area increases and the concentration shrinks. Consequently, a larger concentration difference develops between the constant Dirichlet boundary condition at the top of the domain and the adjacent interior of the domain. This leads to an increased flux of particles into the domain. When the displacement is relaxed, the concentration rises again. After some time (see Figure \ref{fig:Macroscopic_Diffusion} at $t = 25$), the concentration can even increase to values that are higher than the Dirichlet boundary condition $\widehat{c} = 1$. This would not be possible in the case of a non-deforming domain.
\section{Sensitivity of Transport Processes with respect to Changes in the Deformation} 
\label{section:sensitivity_of_transport_process}
\begin{figure}
\centering
\begin{subfigure}{0.49\textwidth}
\includegraphics[width = \textwidth]{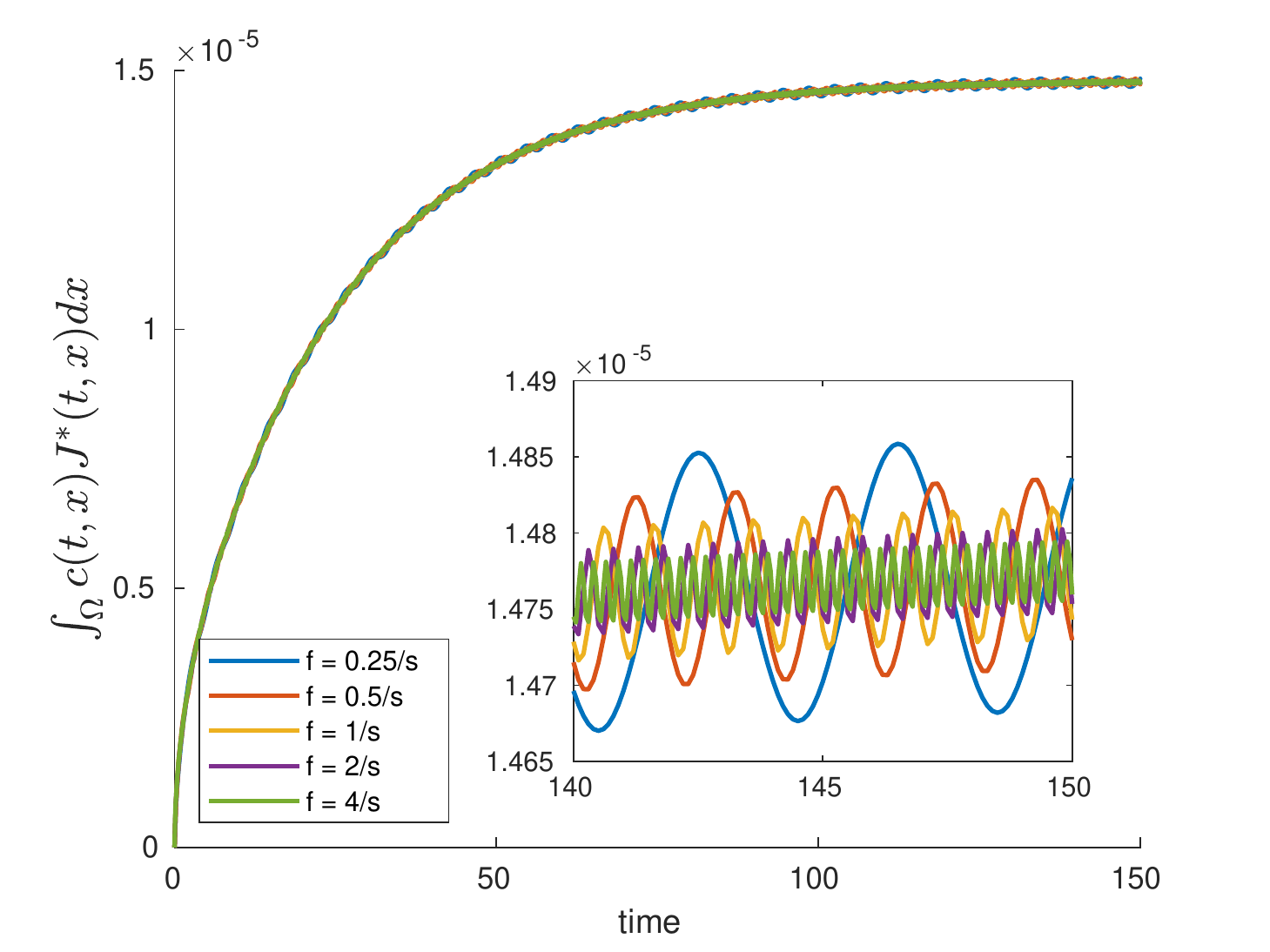}
\caption{frequency dependence}
\label{fig:sensitivity_experiments:a}
\end{subfigure}
\begin{subfigure}{0.49\textwidth}
\includegraphics[width = \textwidth]{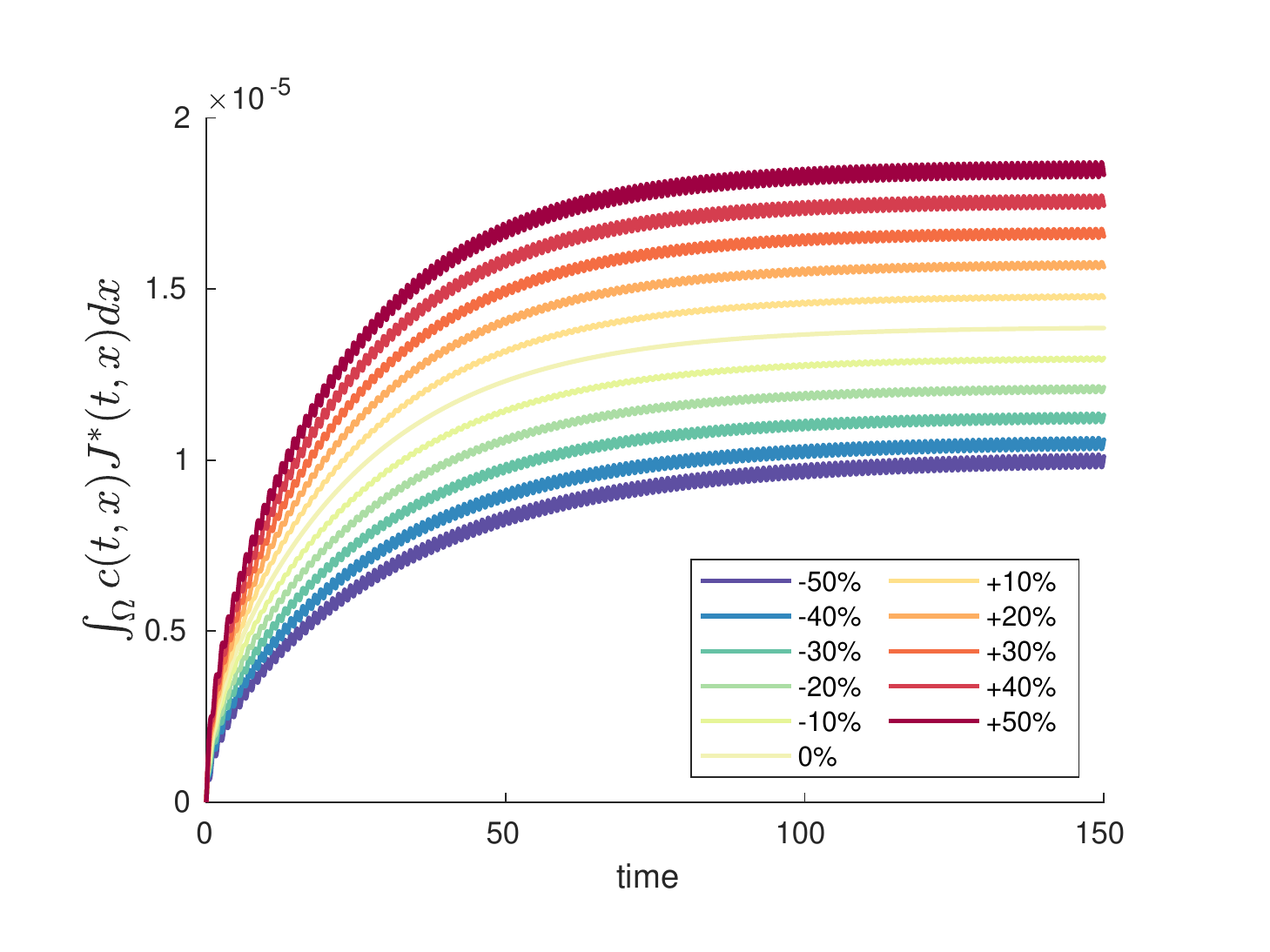}
\caption{amplitude dependence}
\label{fig:sensitivity_experiments:b}
\end{subfigure}
\caption{Sensitivity Experiments}
\label{fig:sensitivity_experiments}
\end{figure}
\begin{figure}
\centering
\begin{subfigure}[b]{0.49\textwidth}
\includegraphics[width = \textwidth]{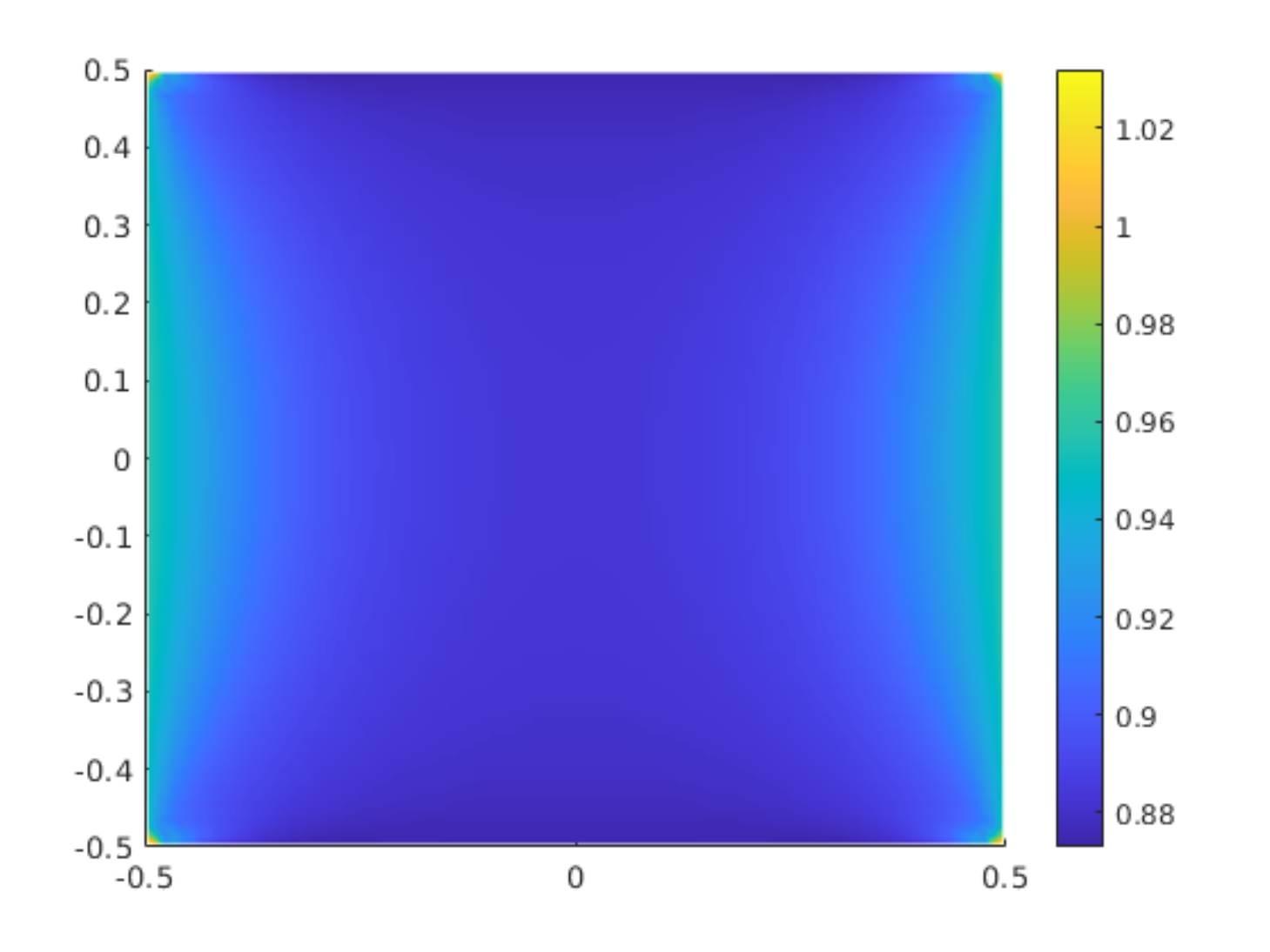}
\caption{$\Jhom(t = 0.5, \cdot)$}
\label{fig:J_star}
\end{subfigure}
\begin{subfigure}[b]{0.49\textwidth}
\includegraphics[width = \textwidth]{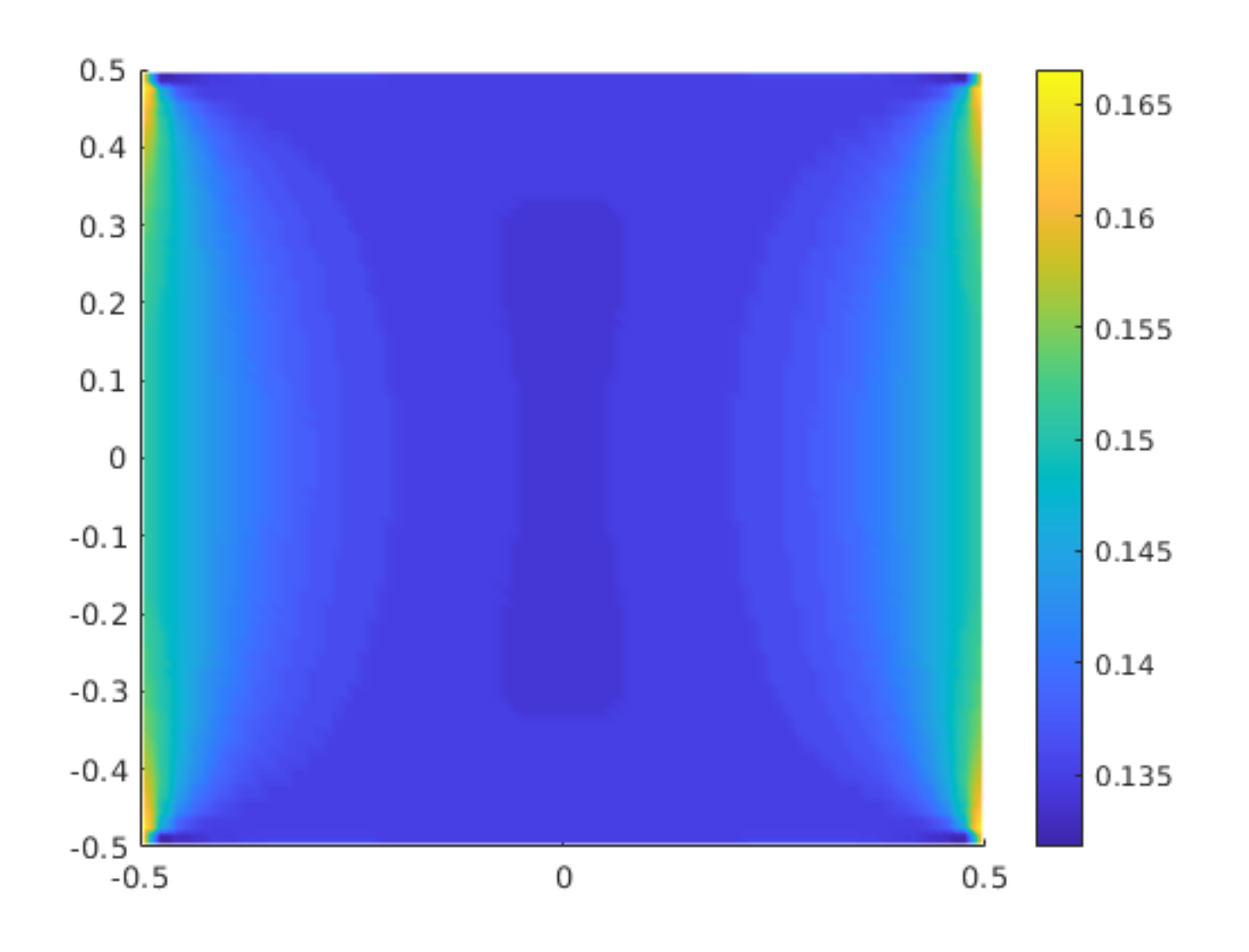}
\caption{$\Dhom_{11}(t = 0.5, \cdot)$}
\label{fig:D_star_11}
\end{subfigure}%
\begin{subfigure}[b]{0.49\textwidth}
\includegraphics[width = \textwidth]{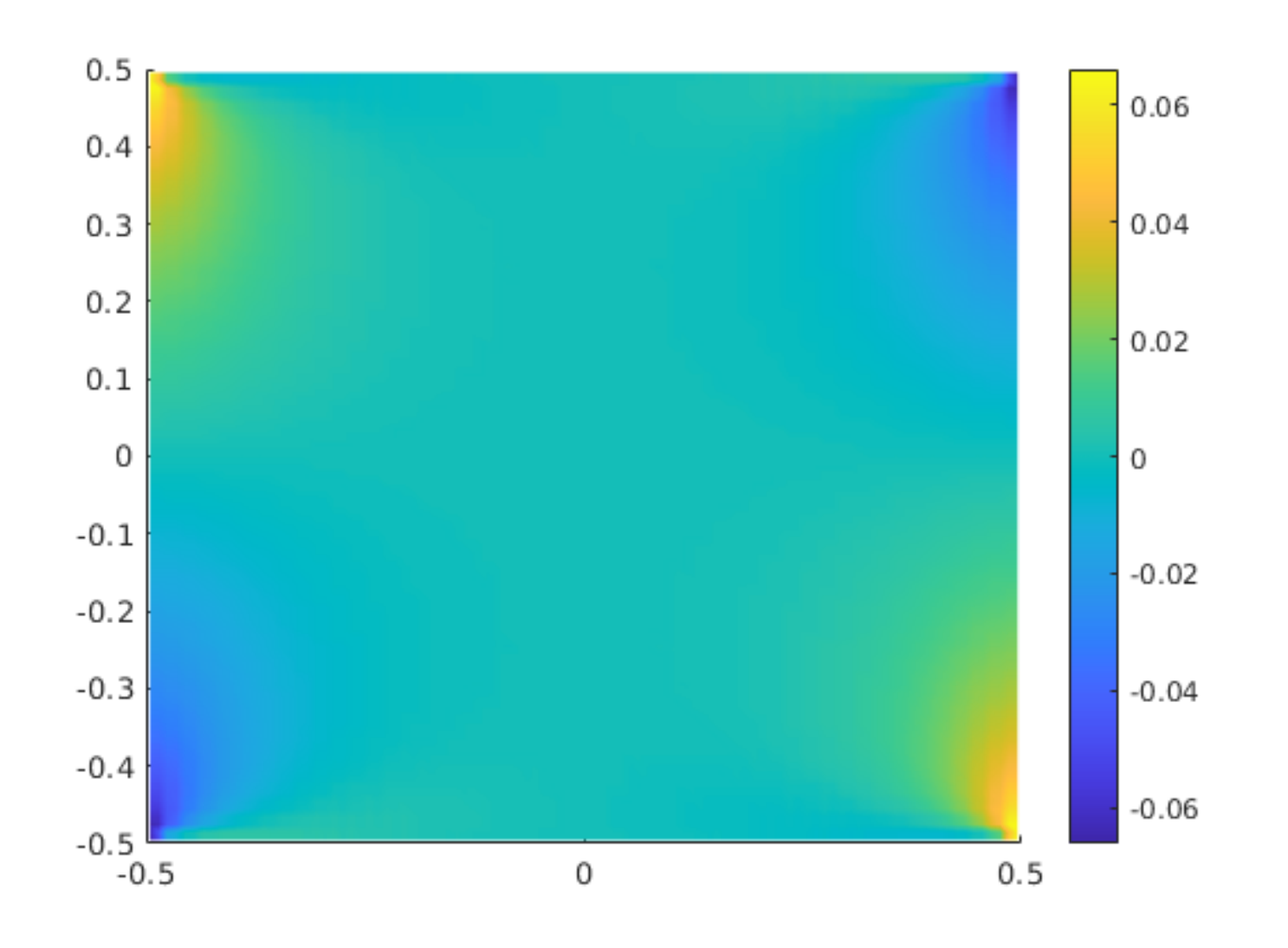}
\caption{$\Dhom_{12}(t = 0.5, \cdot)$}
\label{fig:D_star_12}
\end{subfigure}
\begin{subfigure}[b]{0.49\textwidth}
\includegraphics[width = \textwidth]{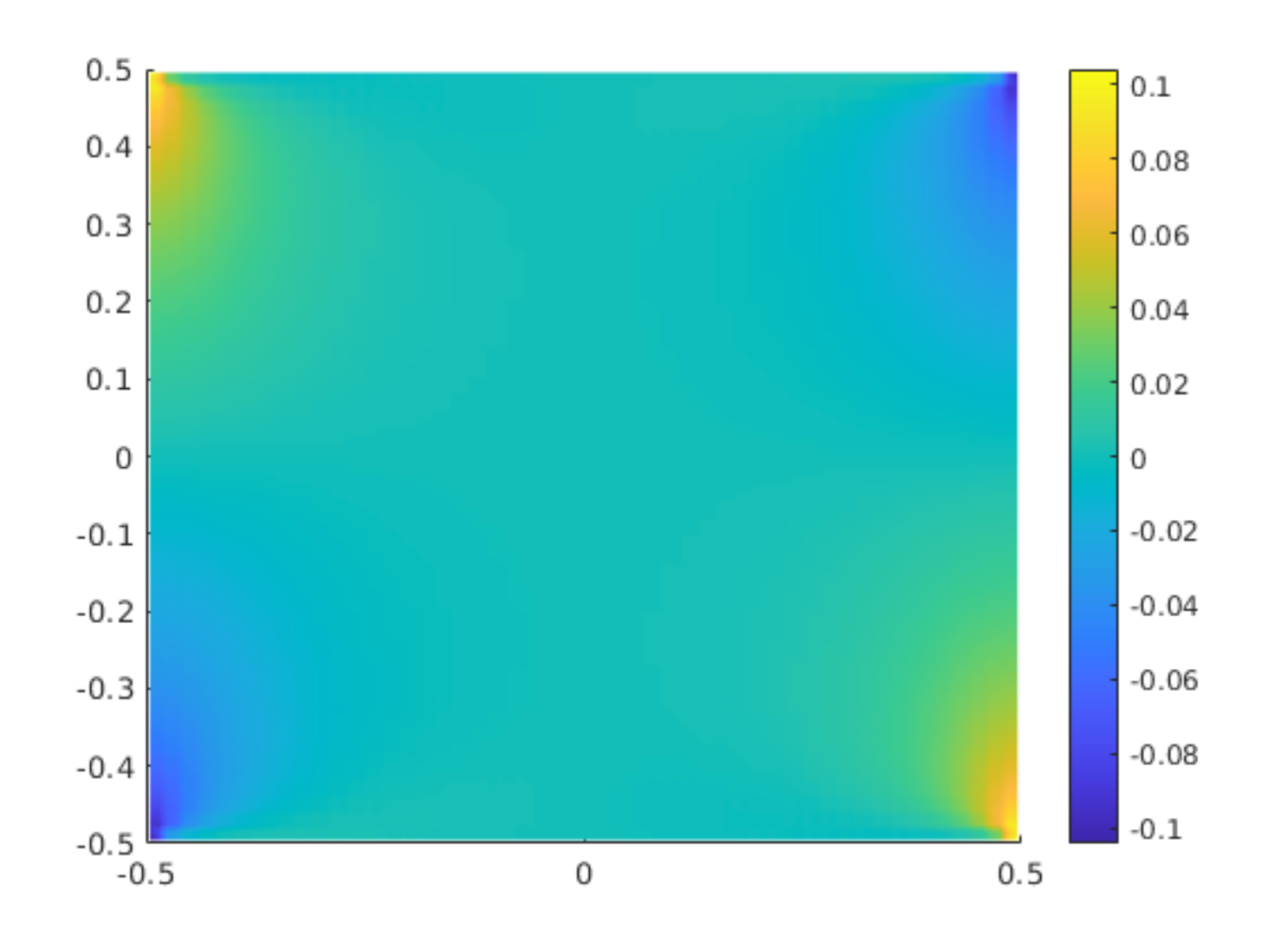}
\caption{$\Dhom_{21}(t = 0.5, \cdot)$}
\label{fig:D_star_21}
\end{subfigure}%
\begin{subfigure}[b]{0.49\textwidth}
\includegraphics[width = \textwidth]{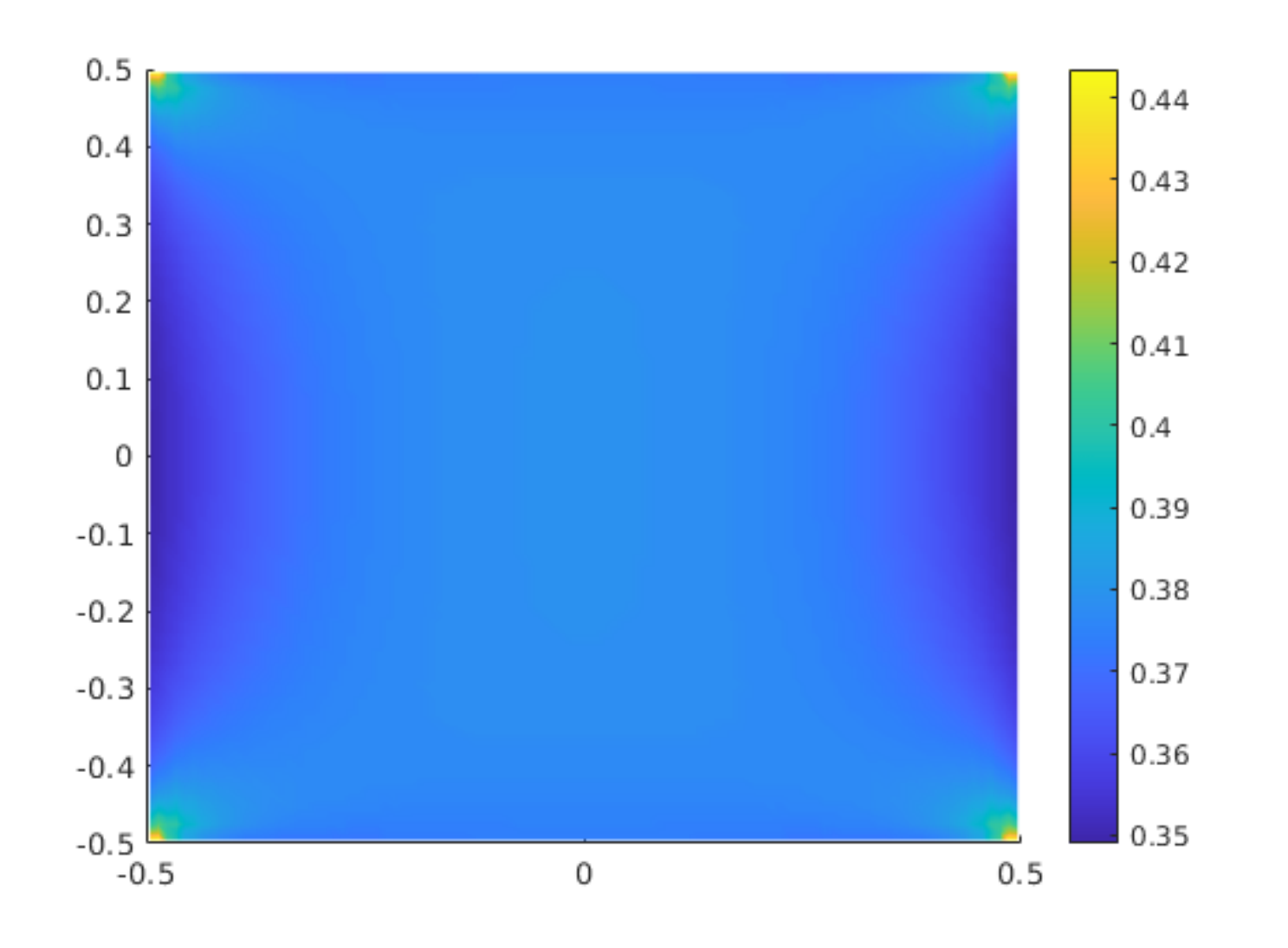}
\caption{$\Dhom_{22}(t = 0.5, \cdot)$}
\label{fig:D_star_22}
\end{subfigure}%
\caption{{\bf Plots of homogenized quantities for maximal lateral displacement} with amplitude $a = 0.25$ at $t = 0.5$. }
\label{fig:effective_coefficients}
\end{figure}
Finally, we analyze the system's sensitivity with respect to experimental parameters which we can modulate by manipulating the boundary conditions of the elasticity problem. We adopt again the setup from the model problem in section \ref{subsection:The_Model_Problem}, i.e. the deformation of the domain is driven by the time-periodic Dirichlet boundary condition $\mathbf{h}$ (cf. eq. (\ref{eq:elasticity_dirichlet_bc_constant_front})), while a diffusing substance spreads into the domain, originating from the Dirichlet boundary condition of the diffusion problem at the upper part of the boundary of $\Omega$.  Since we are concerned with the effect of the deformation on transport processes, we monitor the quantity $M(t) := \int_{\Omega} c(t,x) \Jhom(t,x) \,dx$ over the course of each of the following computer experiments. $M(t)$ can be associated with the mass or the number of particles that accumulates over time in the domain. Different experimental scenarios arise from varying the frequency $f$ or the amplitude $a$ in \eqref{eq:elasticity_dirichlet_bc_constant_front}. By comparing $M(t)$ for different simulations we gain insights on the sensitivity of the transport of the diffusing substance with respect to the deformation parameters. Figure \ref{fig:sensitivity_experiments:a} shows how the mass of the diffusing substance accumulates in the domain for different frequencies of the lateral displacement. The curves indicate no significant qualitative difference at first glance. Only a closer look reveals slight deviations in the small mass fluctuations. The question may arise why the mass in the domain can temporarily decrease during the advanced stage of the simulations. The answer is already partly given in section \ref{subsection:Solution_of_the_Diffusion_Subproblem}: when the domain is released back into the relaxed configuration, the concentration inside the domain may reach values that are greater than the prescribed Dirichlet boundary condition, consequently leading to an outflow of particles across that very boundary. The system clearly exhibits a higher sensitivity with respect to the amplitude of the lateral displacement, as can be seen from Figure \ref{fig:sensitivity_experiments:b}. There, the mass $M(t)$ is plotted for varying amplitudes of the lateral displacement. The percentage in the legend encodes the maximal length of the domain $\Omega^0(t)$ in lateral direction during each experiment in relation to the edge length of $\Omega$. Positive percentages indicate cyclic extension and relaxation into the initial state whereas negative percentages indicate cyclic compression and subsequent relaxation. From the experiments it is obvious that deforming the domain with increasing amplitude leads to the accumulation of more mass in the domain. This effect can be attributed to the increased area of the extended configuration: During expansion, the concentration shrinks, hence leading to a greater concentration difference between the constant Dirichlet boundary condition and the concentration inside the domain, ultimately resulting in an increased flux into the domain across the respective boundary. The converse is true for compression of the domain, leading to the opposite effect.\\
Visualization of the homogenized quantities $\Jhom$ and $\Dhom$ in $\Omega$ for a fixed time point when the physical domain $\Omega^0(t)$ is in the maximally deformed configuration, sheds more light into the properties of the diffusion process during the course of the simulations. In Figures \ref{fig:J_star}--\ref{fig:D_star_22}, $\Jhom$ and the components of $\Dhom$ are plotted for the case where $a = 0.25$, corresponding to 50\% extension from Figure \ref{fig:sensitivity_experiments:b}. Unsurprisingly, the homogenized determinant $\Jhom$, which keeps track of volume changes of the domain, is increased everywhere (note that $\Jhom$ for the relaxed configuration is constant over space with $\Jhom = \lvert Y \rvert = 5/9$ due to the heterogeneity of the underlying microscopic problem). From Figures \ref{fig:D_star_11} and \ref{fig:D_star_22} we deduce that lateral stretching of the domain facilitates diffusion in $x_2$-direction while decelerating diffusion in $x_1$-direction when compared to the values $\Dhom_{11} = \Dhom_{22} = 0.184577$ (cf. Table \ref{tab:comparison_of_initial_and_effective_quantities}) in the non-deformed configuration. Additionally, we note that anisotropy effects are introduced due to non-zero off-diagonal entries $\Dhom_{21}$ and $\Dhom_{21}$. This can be attributed to the deformation of the domain, as we have seen (cf. Table \ref{tab:comparison_of_initial_and_effective_quantities}) that the isotropy property of the initial diffusion coefficient $\hat{\mathbf{D}}$ \eqref{eq:diffusion_coefficient_for_numerics} transfers to the homogenized coefficient, when the domain is not transformed, due to the symmetry of $Y$.
\section{Conclusion}
In this work, an effective model for transport processes in elastically deformable, heterogeneous meda has been derived and investigated. Starting from a microscopic description of an elasticity-diffusion problem, which contained the underlying heterogeneity of the modelled process explicitly  in a mixed \textit{Lagrangian}/\textit{Eulerian} framework, an upscaled model, formulated on a fixed, homogeneous, reference domain on the macroscale, has been derived. The underlying heterogeneity is addressed via cell problems for the deformation and the diffusion process. Although the initial problem is posed with constant coefficient functions, a transformation of the diffusion problem onto the fixed reference domain introduces a diffusion coefficient for the cell problems which depends not only on the microscopic space-variable, but also on time and on the macroscopic space-variable. Solving individual cell problems and computing effective coefficients in each time step and for each quadrature point for the effective diffusion problem poses major difficulties for the numerical simulation of the upscaled problem. There is some potential for optimization by computing cell solutions in parallel or by solving fewer, representative cell problems which a chosen by an adaptive algorithm. \\
Investigations associated to the numerical solutions of the upscaled problem such as convergence tests and interpretation of solutions within the context of realistic experimental setups, inspired by applications, are used to verify the implemented model and give insight on the qualitative and quantitative behavior of the model. The simulations emphasize, that the deformation can significantly influence the transport processes and therefore shows the necessity to consider e.g.  the effect of respiratory movement when investigating transport processes on the cellular level in the lung or related problems. \\
The present model has to be understood as a starting point for further research in different directions. As a first generalization, it is reasonable to expand the model by including also an influence of the diffusing substance on the elastic properties of the solid, leading to a fully coupled system. In the medical context, this may be especially interesting when the effect of fibrosis on the functionality of the lung has to be considered. Furthermore, until now, the effects of the perforations of the microscopic domain are only considered in the sense that they introduce the heterogeneity of the problem. From a biological perspective, the space that is not occupied by cells can be interpreted as the extracellular matrix, which, as a first approximation, can be modelled by a fluid. Starting from this, a more comprehensive model could be derived, governing also the effects of fluid-structure-interaction. Finally, as we state in the introduction, the model is inspired by the lung and a related biomimetic microdevice that shares common features with the lung. Although the structures of these examples are based on thin layers (e.g. air-blood-barrier, consisting only of two cellular layers and the basal lamina), we chose to formulate the model for a domain that spreads in all space dimensions for the sake of simplicity. Nonetheless, for a model governing an elasticity-diffusion process for a thin layer, an additional limit process could be introduces by sending the thickness of the domain in one direction to zero.

\bmhead{Acknowledgments} This research is supported by SCIDATOS (Scientific Computing for Improved Detection and Therapy of Sepsis). SCIDATOS is a collaborative project funded by the Klaus Tschira Foundation, Germany (Grant Number 00.0277.2015) and provided in particular the funding for the research of the first and the second author. 

%
\bibliography{bibliography}


\end{document}